\documentclass[12pt]{amsart}
\usepackage{amssymb,latexsym,amstext,amsmath,epsfig,ifthen}

\numberwithin{equation}{section}

\newtheorem{theorem}{Theorem}[section]
\newtheorem{proposition}[theorem]{Proposition}
\newtheorem{corollary}[theorem]{Corollary}
\newtheorem{lemma}[theorem]{Lemma}
\newtheorem{conjecture}[theorem]{Conjecture}
\newtheorem{problem}[theorem]{Problem}

\theoremstyle{definition}

\newtheorem{remark}[theorem]{Remark}
\newtheorem{example}[theorem]{Example}
\newtheorem{definition}[theorem]{Definition}

\def\RR{\mathbb{R}}
\def\ZZ{\mathbb{Z}}
\def\CC{\mathbb{C}}

\def\Upper{\mathcal{U}}

\def\AA{\mathcal{A}}

\def\FFcal{\mathcal{F}}

\def\ZZ{\mathbb{Z}}
\def\CC{\mathbb{C}}
\def\RR{\mathbb{R}}
\def\QQ{\mathbb{Q}}

\def\cc{\mathbf{c}}
\def\ii{\mathbf{i}}

\def\xx{\mathbf{x}}
\def\ex{\mathbf{ex}}

\newcommand{\bmat}[4]{\left[\!\!\begin{array}{cc}
#1 & #2 \\ #3 & #4 \\
\end{array}\!\!\right]}

\newcounter{colorversion}

\usepackage{color}
\newcommand{\verylight}[1]{\color{lightgray}{#1}\color{black}}
\definecolor{light}{rgb}{1,0.1,0.1}        
\definecolor{dark}{cmyk}{1,0.4,0,0}  
\definecolor{bulletcolor}{cmyk}{1,0.3,1,0}  
\definecolor{lightgray}{rgb}{0.8,0.8,0.8}
\definecolor{gray}{rgb}{0.7,0.7,0.7}

\ifthenelse{\thecolorversion=1}
{

}{}

\begin{document}

\ \vspace{-.1in}

\title[Cluster algebras: CDM-03 Notes]{{\ }\\Cluster algebras:\\[.1in]
Notes for the CDM-03 conference}

\author{Sergey Fomin}
\address{Department of Mathematics, University of Michigan,
Ann Arbor, MI 48109, USA} \email{fomin@umich.edu}

\author{Andrei Zelevinsky}
\address{\noindent Department of Mathematics, Northeastern University,
 Boston, MA 02115, USA}
\email{andrei@neu.edu}



\date{September 21, 2004}

 \thanks{Research 
supported 
by NSF (DMS) grants 0245385 (S.F.) and 0200299~(A.Z.).}

\begin{abstract}
This is an expanded version of the notes of our lectures
given at the conference \textsl{Current Developments in Mathematics~2003}
held at Harvard University on November 21--22, 2003.
We present an overview of the main definitions, results and
applications of the theory of cluster algebras.
\end{abstract}

\maketitle


\tableofcontents

\newpage

\section{Introduction}


Cluster algebras, first introduced and studied in~\cite{ca1,ca2,ca3},
are a class of axiomatically defined commutative
rings equipped with a distinguished set of generators (cluster
variables) grouped into overlapping subsets (clusters) of the same
finite cardinality.
The original motivation for this theory lied in the desire to create an
algebraic framework for total positivity
and canonical bases in semisimple algebraic groups.

Since its inception, the theory of cluster algebras has developed
several interesting connections and applications:

\begin{itemize}
\item
Discrete dynamical systems based on rational
recurrences~\cite{fz-laurent, carroll-speyer, speyer}.
\item
$Y$-systems in thermodynamic
Bethe Ansatz~\cite{yga}.
\item
Generalized associahedra
associated with finite root systems~\cite{yga,cfz}.
\item
Quiver representations~\cite{bmrrt, bmr, ccs, mrz}.
\item
Grassmannians, projective configurations and
their tropical analogues~\cite{scott,speyer-williams}.
\item
Quantum cluster algebras, Poisson geometry and Teichm\"uller theory~\cite
{fock-goncharov1, fock-goncharov2, gsv1, gsv2, qca}.
\end{itemize}

In these lectures, we concentrate on the following aspects.
Sections~\ref{sec:tp-dbc} and~\ref{sec:triangulations}
set the groundwork for the future theory by supplying
a family of motivating examples.
Specifically,
Section~\ref{sec:tp-dbc}\footnote{
This section contains material that was not covered in CDM-2003
lectures due to time constraints.
We feel, however, that it is important to
present in these notes what was historically the main motivating example for the
creation of cluster algebra theory.
}
is devoted to total positivity and geometry of double Bruhat cells
in semisimple groups (a closely related connection to canonical bases
was discussed in~\cite{camblec}),
while a more elementary setup of Section~\ref{sec:triangulations}
involves Ptolemy relations, combinatorics of triangulations, and
a Grassmannian of $2$-dimensional subspaces.

Section~\ref{sec:fundam} introduces cluster algebras in earnest,
and discusses some basic results, open problems and conjectures.
Section~\ref{sec:finite-type} focuses on cluster algebras of finite
type, including their complete classification and their combinatorics,
which is governed by generalized associahedra,
a family of convex polytopes associated with finite crystallographic
root systems.

We present no proofs, each time referring the reader to primary
sources.
Some of the material in these notes is based on the coverage in the
earlier surveys~\cite{pcmi,zel-korea}.
Pictures were borrowed from~\cite{yga,cfz,pcmi}.


\section{Total positivity and double Bruhat cells}
\label{sec:tp-dbc}

\subsection{Total positivity in semisimple Lie groups}

A matrix is \emph{totally positive} \ (resp., \emph{totally nonnegative})
if all its minors are positive (resp., nonnegative) real numbers.
The first systematic study of these classes of matrices was undertaken
in the 1930s by I.~J.~Schoenberg, and by F.~R.~Gantmacher and
M.~G.~Krein.
References to their papers and a discussion of further connections and
applications of totally positive matrices (with an emphasis on
algebraic and combinatorial aspects) can be found in~\cite{fz-intel}.
The interest in the subject intensified in the last decade
due in large part to the discovery by G.~Lusztig~\cite{lusztig}
of a surprising connection between total positivity and
canonical bases for quantum groups.
Among other things, Lusztig extended the subject~by defining
the \emph{totally positive variety} $G_{> 0}$ and the
\emph{totally nonnegative variety} $G_{\geq 0}$
inside every complex reductive group~$G$.
These ideas were further developed
in~\cite{bfz-tp, bz-schubert, fz-double, fz-osc}.
We are going to present some of the results obtained in those papers.
For technical reasons, we restrict ourselves to the case where~$G$
is semisimple and simply connected.
Classical total positivity theory is recovered for
$G = SL_{r+1}(\CC)$, the group of complex
$(r+1) \times (r+1)$ matrices with determinant~$1$.

Lusztig defined $G_{> 0}$ and $G_{\geq 0}$ parametrically.
As shown in~\cite{fz-osc}, Lusztig's original definition is
equivalent to setting
\begin{align}
\label{eq:gtp-gtnn}
G_{> 0}& = \{x \in G: \Delta_{\gamma,\delta}(x) > 0
\,\, \text{for all} \,\, \gamma,\delta\}, \\
\nonumber
G_{\geq 0} &= \{x \in G: \Delta_{\gamma,\delta}(x) \geq 0
\,\, \text{for all} \,\, \gamma,\delta\},
\end{align}
for a family of regular functions~$\Delta_{\gamma,\delta}$
on~$G$ called \emph{generalized minors}.
To be more specific, let $r$ be the rank of~$G$,
let $\omega_1, \dots, \omega_r$ be the
fundamental weights, and let~$W$ be the Weyl group.
(The standard terminology and notation used here is explained in more detail
in~\cite{fz-double}.)
The indices~$\gamma$ and~$\delta$ of a generalized minor are two
elements of the $W$-orbit of the same fundamental weight~$\omega_i$, i.e., two
extremal weights in the same fundamental representation~$V_{\omega_i}$
of~$G$.
The corresponding minor $\Delta_{\gamma,\delta}$ is a suitably
normalized matrix element of~$V_{\omega_i}$ associated with
the weights~$\gamma$ and~$\delta$.

In the case $G = SL_{r+1}(\CC)$,
the regular functions  $\Delta_{\gamma,\delta}$ specialize to the ordinary minors
(i.e., determinants of square submatrices) as follows.
The Weyl group $W$ is identified
with the symmetric group~$S_{r+1}$, and
$V_{\omega_i} = \bigwedge^i \CC^{r+1}$, the $i$th exterior power
of the standard representation.
All the weights of $V_{\omega_i}$
are extremal, and are in bijection with the $i$-subsets of
$[1,r+1] = \{1,\dots,r+1\}$, so that~$W = S_{r+1}$
acts on them in a natural way, and $\omega_i$ corresponds to~$[1,i]$.
If~$\gamma$ and~$\delta$ correspond to $i$-subsets $I$ and $J$, respectively,
then $\Delta_{\gamma,\delta} = \Delta_{I,J}$ is the minor with
the row set~$I$ and the column set~$J$.

The natural geometric framework for studying the varieties
$G_{> 0}$ and $G_{\geq 0}$ is provided by the
\emph{double Bruhat cells}
\begin{equation}
\label{eq:double-cell}
G^{u,v} = B u B \cap B_- v B_- \,;
\end{equation}
here $u,v \in W$, and~$B$ and~$B_-$ are two opposite Borel
subgroups in~$G$.
(For $G = SL_{r+1}(\CC)$, the standard choice for~$B$
(resp.,~$B_-$) is the subgroup of upper- (resp., lower-) triangular
matrices.)
The group~$G$ has two \emph{Bruhat decompositions}, with respect
to~$B$ and $B_-\,$:
\begin{equation}
G = \bigcup_{u \in W} B u B = \bigcup_{v \in W} B_- v B_- \ .
\end{equation}
Thus, $G$ is the disjoint union of all double Bruhat
cells\footnote{The term ``cell" is somewhat misleading since the
  topology of $G^{u,v}$ may be non-trivial.}.

Recall that the Weyl group~$W$ is a Coxeter group generated by \emph{simple
reflections} $s_1, \dots, s_r$.
For $w \in W$, the \emph{length} $\ell = \ell(w)$ is the shortest length of a
sequence of indices $\ii = (i_1, \dots, i_\ell)$ such that
$w = s_{i_1} \cdots s_{i_\ell}$.
Such a sequence is called a \emph{reduced word} for~$w$.
For $G = SL_{r+1}(\CC)$, the simple reflection
$s_i \in W = S_{r+1}$ is the transposition of~$i$ and~$i+1$,
and $\ell(w)$ is the number of inversions of a permutation~$w$.

\begin{theorem}[{\cite[Theorem~1.1]{fz-double}}]
The variety $G^{u,v}$ is isomorphic to a Zariski open subset of an
affine space of dimension
$r+\ell(u)+\ell(v)$. 
\end{theorem}

\begin{example}
Let $w_\circ$ be the longest element in~$W$.
Then $G^{w_\circ, w_\circ}$ is the open double Bruhat cell given by:
\begin{equation}
G^{w_\circ, w_\circ} = \left\{x \in G:
\begin{array}{l}
\Delta_{\omega_i, w_\circ \omega_i}(x) \neq 0 \\
\Delta_{w_\circ \omega_i, \omega_i}(x) \neq 0
\end{array}
\text{for all}
\,\, i \in [1,r]\right\}.
\end{equation}
\end{example}

Following \cite{fz-double}, we define the \emph{totally positive part}
of $G^{u,v}$ by setting
\begin{equation}
G^{u,v}_{> 0} = G^{u,v} \cap G_{\geq 0} \ .
\end{equation}

\begin{proposition}
\label{pr:tp-part-biggest-cell}
The totally positive part of the open double Bruhat cell is the
totally positive variety:
$G^{w_\circ, w_\circ}_{> 0} = G_{> 0}\,$.
\end{proposition}

\begin{example}
As a nice exercise in linear algebra, the reader is invited to
check Proposition~\ref{pr:tp-part-biggest-cell} for $G = SL_3(\CC)$.
In other words: if all minors of a matrix $x=(x_{ij}) \in SL_3(\CC)$
are real and nonnegative, and the minors
$$\Delta_{1,3}(x) = x_{13}, \,\, \Delta_{12,23}(x), \,\,
\Delta_{3,1}(x) = x_{31}, \,\, \Delta_{23,12}(x)$$
are nonzero, then all the minors of~$x$ are positive.
\end{example}

It turns out that the variety~$G^{u,v}_{> 0}$ is isomorphic to
$\RR_{>0}^{r+\ell(u)+\ell(v)}$.
This was first demonstrated by G.~Lusztig~\cite{lusztig} by
introducing local coordinates in~$G^{u,v}$ which consist of rational
functions that are not necessarily regular.
To sharpen this result, we will need the following notion,
which is a slight variation of the definition given in~\cite{fz-double}.

\begin{definition}
\label{def:tp-basis}
A \emph{TP-basis} for $G^{u,v}$ is a collection of regular
functions $F= \{f_1, \dots, f_m\} \subset \CC[G^{u,v}]$ with the
following properties:
\begin{enumerate}
\item[(i)]
The functions $f_1,\dots,f_m$ are algebraically independent and
  generate the field of rational functions $\CC(G^{u,v})$; in
  particular, \linebreak
$m = r+\ell(u)+\ell(v)$.
\item[(ii)]
The map $(f_1, \dots, f_m): G^{u,v} \to \CC^m$ restricts to a
  biregular isomorphism $U(F) \to (\CC_{\neq 0})^m$, where
\begin{equation}
\label{eq:U(F)}
U(F)=\{x \in G^{u,v} : \text{$f_k(x) \neq 0$ for all $k \in [1,m]$}\}.
\end{equation}
\item[(iii)]
The map $(f_1, \dots, f_m): G^{u,v} \to \CC^m$
  restricts to an isomorphism $G^{u,v}_{> 0} \to \RR_{> 0}^m$.
\end{enumerate}
\end{definition}

Condition (iii) can be rephrased as saying that 
$f_1,\dots,f_m$ provide a \emph{total positivity criterion}
in~$G^{u,v}$: an element $x\in G^{u,v}$ is totally nonnegative
(i.e., $x\in G_{\geq 0}$) if and only if $f_k(x)>0$ for
all~$k\in[1,m]$.

The following theorem from~\cite{fz-double} extends the
results previously obtained in~\cite{bfz-tp, bz-schubert}.

\begin{theorem}
\label{th:tp-basis}
Each reduced word~$\ii$ for $(u,v) \in W \times W$ gives rise to a
TP-basis $F_\ii$ for~$G^{u,v}$ consisting of generalized minors:
\begin{equation}
\label{eq:Fi}
F_\ii =
\{\Delta_{\gamma_k, \delta_k}: k \in [1,m]\};
\end{equation}
here $m = r + \ell(u) + \ell(v)$,
and the weights $\gamma_k, \delta_k$ are defined explicitly in terms of
$\ii$ in Definition~\ref{def:gamma-delta-k} below.
\end{theorem}

\begin{definition}
\label{def:gamma-delta-k}
We will represent~$\ii$ as a sequence of indices $(i_1, \dots, i_m)$
from the set
\[
-[1,r] \cup [1,r]
\stackrel{\rm def}{=}
\{-1, \dots, -r\}\cup\{1,\dots,r\}
\]
such that
$i_j = j$ for $j \in [1,r]$, and
\begin{equation}
\label{eq:ii-products}
s_{-i_{r+1}} \cdots s_{-i_m} = u, \quad
s_{i_{r+1}} \cdots s_{i_m} = v \, ,
\end{equation}
with the convention
$s_{-i} = 1$ for $i \in [1,r]$.
In this notation, we define:
\begin{equation}
\label{eq:gamma-k}
\gamma_k = s_{-i_1} \cdots s_{-i_k} \omega_{|i_k|}, \quad
\delta_k = s_{i_m} \cdots s_{i_{k+1}} \omega_{|i_k|} \, .
\end{equation}
\end{definition}

\pagebreak[3]

\begin{example}
\label{example:SL3-w0w0}
Let $G=SL_3(\CC)$, and let $u=v=w_\circ=s_1s_2s_1=s_2s_1s_2$ be the order-reversing
permutation (the element of maximal length in the symmetric group~$W\!=\!S_3$).
Take $\ii=(1,2,1,2,1,-1,-2,-1)$.
Then
\begin{equation}
\label{eq:f-sl3}
F_\ii = \{\Delta_{1,3}, \Delta_{12,23}, \Delta_{1,2},
\Delta_{12,12}, \Delta_{1,1}, \Delta_{2,1}, \Delta_{23,12},
\Delta_{3,1}\}.
\end{equation}
(The minors on the right-hand side are listed in the natural order,
i.e., $f_1=\Delta_{1,3}\,$, \dots, $f_8=\Delta_{3,1}\,$.)
The reader is invited to check that these $8$ minors provide a total
positivity criterion in $SL_3$:
a $3\times 3$ matrix $x$ with determinant~$1$ is totally positive if and
only if the evaluations of these $8$ minors at $x$ are positive.
\end{example}

\subsection{
Coordinate rings of double Bruhat cells}
\label{sec:exchanges-double-cells}

Theorem~\ref{th:tp-basis} was obtained a couple of years before the
discovery of cluster algebras. In retrospect, it can be seen as the
first step towards introducing the (upper) cluster algebra structure in the
coordinate ring~$\CC[G^{u,v}]$.
(Each collection $F_\ii$ would give rise to a
cluster in the future theory.)
An impetus for introducing this structure came from a realization that
working with the families~$F_\ii$ may not be sufficient for answering
some natural questions about double Bruhat cells.
Here is one such question: describe and enumerate the connected
components of the real part of~$G^{u,v}$.
An important special case of this problem, with $u = 1$ and
$v = w_\circ$ for $G=SL_{r+1}(\CC)$, was solved in~\cite{ssv1,ssv2};
see also related work~\cite{rie1,rie2}.
The general case was handled in \cite{z-imrn} using results
and ideas from \cite{ssvz} and the earlier papers mentioned above.
(For follow-ups see~\cite{gsv0,seven}.)
The solution in~\cite{z-imrn} utilized the following general approach,
which goes back to~\cite{ssv1}: try to find
a ``simple" Zariski open subvariety~$U\subset G^{u,v}$ such that
the codimension in $G^{u,v}$ of the complement of~$U$ is greater
than~$1$.
If furthermore $U$ is ``compatible" with the real part
of $G^{u,v}$, then replacing  $G^{u,v}$ by~$U$ does not change
the structure of connected components of the real part.
In many cases (including the one covered in~\cite{ssv1,ssv2}) one
can take
\begin{equation}
\label{eq:U=UU(F)}
U=\bigcup_\ii U(F_\ii),
\end{equation}
where $\ii$ runs over all reduced words for~$(u,v)$,
and $U(F_\ii)$ is defined by~\eqref{eq:U(F)}.
In general, however, the complement of the open subvariety
\eqref{eq:U=UU(F)} has codimension~$1$, so one needs something else.

Here is the main construction from~\cite{z-imrn}.
Fix a pair~$(u,v) \in W \times W$ and a reduced word
$\ii$ for~$(u,v)$, written in the form $\ii = (i_1, \dots, i_m)$ as in
Definition~\ref{def:gamma-delta-k}.
For $k \in [1,m]$, abbreviate~$f_k = \Delta_{\gamma_k, \delta_k}$,
where $\gamma_k$ and $\delta_k$ are given by \eqref{eq:gamma-k}.
We say that an index~$k \in [1,m]$ is \emph{$\ii$-exchangeable} if
$r < k \leq m$, and $|i_p| = |i_k|$ for some $p > k$.
Denote by $\ex = \ex_\ii \subset [1,m]$ the subset of
$\ii$-exchangeable indices.

\pagebreak[3]

\begin{lemma}
The subset
\begin{equation}
\label{eq:cc-bdc}
\cc=\{f_k: k \in [1,m] - \ex\} \, \subset \, F_\ii
\end{equation}
depends only on~$u$ and~$v$, not on the particular choice of a
reduced word~$\ii$.
In particular, the cardinality~$n$ of~$\ex$ depends only on~$u$ and~$v$.
\end{lemma}

The set $\cc$ can be described explicitly in terms of $u$ and~$v$ as
follows:  $f_j = \Delta_{\omega_j, v^{-1}\omega_j}$
for $j \in [1,r]$, and $f_k = \Delta_{u \omega_i,\omega_i}$
if~$i_k$ is the last occurrence of $\pm i$ in~$\ii$.

\begin{example}
\label{example:SL3-w0w0-2}
In the special case of Example~\ref{example:SL3-w0w0}
(that is,  $G=SL_3(\CC)$ and $u=v=w_\circ$),
we have $\ex = \{3,4,5,6\}$ and
\[
\cc=\{f_1,f_2,f_7,f_8\}=
\{\Delta_{1,3}, \Delta_{12,23}, \Delta_{23,12}, \Delta_{3,1}\}.
\]
\end{example}

\begin{lemma}
\label{lem:cc-nonzero}
Each $f_k\in\cc$  vanishes nowhere in~$G^{u,v}$.
\end{lemma}

By Lemma~\ref{lem:cc-nonzero}, the coordinate ring $\CC[G^{u,v}]$
contains the Laurent polynomial ring~$\CC[\cc^{\pm 1}]$. Our
efforts to understand $\CC[G^{u,v}]$ should therefore concentrate
on the functions~$f_k\notin\cc$, that is, those with exchangeable
indices~$k$.

\begin{theorem}
\label{th:adjacent-clusters}
There is an integer $m \times n$ matrix $\tilde B = \tilde B(\ii) =
(b_{ik})$
(defined explicitly by \eqref{eq:tildeB-entries} below)
with rows labeled by~$[1,m]$ and columns labeled by~$\ex$, such that:
\begin{itemize}
\item
for every~$k \in \ex$, the function
\begin{equation}
\label{eq:exchange-f}
f'_k = \frac
{\prod\limits_{b_{ik}>0} f_i^{b_{ik}}+
\prod\limits_{b_{ik}<0} f_i^{-b_{ik}}}{f_k}
\end{equation}
is regular on~$G^{u,v}$ (i.e., $f_k'\in\CC[G^{u,v}]$),
and the collection
\begin{equation}
\label{eq:Fik}
F_{\ii;k} = F_\ii - \{f_k\} \cup \{f'_k\} \, \subset \, \CC[G^{u,v}]
\end{equation}
is a TP-basis for $G^{u,v}$;
\item
the complement in $G^{u,v}$ of the Zariski open subset $U$ given by
\[
U = U(F_\ii) \, \cup \, \bigcup_{k \in \ex} U(F_{\ii;k})\subset
G^{u,v}
\]
has codimension greater than~$1$.
\end{itemize}
Consequently, $\CC[G^{u,v}] = \CC[U]$.
\end{theorem}


\pagebreak[3]

For the sake of completeness, we provide the definition of the
matrix~$\tilde B = \tilde B(\ii)$ in
Theorem~\ref{th:adjacent-clusters} 
(cf.~\cite[(8.7)]{qca}).
For $p \in [1,m]$ and $k \in \ex$, we have
\begin{equation}
\label{eq:tildeB-entries}
b_{pk} =
\begin{cases}
-\varepsilon (i_k) & \text{if $p = k^-$;} \\
-\varepsilon (i_k) a_{|i_p|, |i_k|} &
\text{if $p < k < p^+ < k^+, \, \varepsilon (i_k) = \varepsilon(i_{p^+})$,}\\
 & \text{or $p < k < k^+ < p^+, \, \varepsilon (i_k) = -\varepsilon (i_{k^+})$;} \\
\varepsilon (i_p) a_{|i_p|, |i_k|} &
\text{if $k < p < k^+ < p^+, \, \varepsilon (i_p) = \varepsilon(i_{k^+})$,}\\
& \text{or $k < p < p^+ < k^+, \, \varepsilon (i_p) = -\varepsilon (i_{p^+})$;} \\
\varepsilon (i_p) & \text{if $p = k^+$;} \\
0 & \text{otherwise,}
\end{cases}
\end{equation}
where we use the following notation and conventions:
\begin{itemize}
\item
$\varepsilon(i) = \pm 1$ for $i \in \pm [1,r]$;
\item
$A = (a_{ij})_{i,j \in [1,r]}$ is the Cartan matrix of~$G$,
that is, the transition matrix between the simple roots and
the fundamental weights:
\begin{equation}
\label{eq:cartan-transition}
\alpha_j = \sum_{i=1}^r a_{ij} \omega_i \ .
\end{equation}
\item
for $k \in [1,m]$, we denote by $k^+ = k^+_\ii$
the smallest index $\ell$ such that $k < \ell \leq m$ and
$|i_\ell| = |i_k|$; if $|i_k| \neq |i_\ell|$ for $k < \ell \leq m$, then
we set $k^+ = m+1$;
\item
$k^- = k^-_\ii$ denotes the index $\ell$ such that $\ell^+=k$;
if such an $\ell$ does not exist, we set $k^- = 0$.
\end{itemize}

\begin{example}
\label{ex:sl3-exchanges}
Let~$G$ and $\ii$ be as in
Examples~\ref{example:SL3-w0w0} and~\ref{example:SL3-w0w0-2}.
Then
\[
\tilde B
= \tilde B(\ii)
=  \begin{pmatrix}
-1 &  0 &  0 &  0\\
1 &  -1 &  0 &  0\\
0 & 1 &  -1 &  0\\
-1 &  0 & 1 &  -1\\
1 &  -1 &  0 & 1\\
0 & 1 &  -1 &  0\\
0 &  -1 &  0 & 1\\
0 &  0 &  0 &  -1
  \end{pmatrix},
\]
where the columns are indexed by $\ex = \{3,4,5,6\}$.
(Cf.\ \cite[Example~3.2]{qca} and \cite[Example~2.5]{ca3}.)
Applying \eqref{eq:exchange-f} to the $f_i$'s in
\eqref{eq:f-sl3}, we get
\begin{align*}
f'_3 &
= \frac{f_2 f_5 + f_1 f_4}{f_3}
= \frac{\Delta_{12,23} \Delta_{1,1} + \Delta_{1,3}\Delta_{12,12}}{\Delta_{1,2}}
= \Delta_{12,23} \, ,\\
f'_4 &
= \frac{f_3 f_6 + f_2 f_5 f_7}{f_4}
= \frac{\Delta_{1,2}\Delta_{2,1}+\Delta_{12,23}\Delta_{1,1}\Delta_{23,12}}{\Delta_{12,12}}
= \Delta_{1,1}\Delta_{23,23} - 1 \, ,\\
f'_5 &= \frac{f_4 + f_3 f_6}{f_5} = \Delta_{22} \, ,\\
f'_6 &= \frac{f_5 f_7 + f_4 f_8}{f_6} =
\Delta_{13,12} \, .
\end{align*}
We note that, while each of the families $F_{\ii;3}, F_{\ii;5}$,
and~$F_{\ii;6}$ is itself of the form $F_{\ii'}$ for some reduced
word~$\ii'$, this is not the case for $F_{\ii;4}$, \linebreak[3]
if only for
the reason that~$f'_4$ is not a minor.
\end{example}

The last statement in Theorem~\ref{th:adjacent-clusters}
implies the following description of the coordinate ring~$\CC[G^{u,v}]$.

\begin{corollary}
\label{cor:CGuv}
The subalgebra $\CC[G^{u,v}]$ of the field of rational functions
$\CC(G^{u,v})$ is the intersection of~$n+1$ Laurent polynomial rings:
\begin{equation}
\label{eq:CGuv-Laurent}
\CC[G^{u,v}] = \CC[F_\ii^{\pm 1}] \cap \bigcap_{k \in \ex}
\CC[F_{\ii;k}^{\pm 1}],
\end{equation}
where the collections~$F_{\ii;k}$ are given by \eqref{eq:Fik}.
\end{corollary}

As we will see in Section~\ref{sec:fundam}, Corollary~\ref{cor:CGuv}
prepares the ground for the concept of an \emph{upper cluster
algebra}.


\section{Ptolemy relations and triangulations
}
\label{sec:triangulations}

In this section, we present another important motivating example
for the concept of cluster algebras:
a family of classically studied commutative algebras~$\AA_n$
($n\in\ZZ_{>0}$).
In the forthcoming classification of cluster algebras of finite type
(see Section~\ref{sec:finite-type}),
$\AA_n$ will turn out to be ``of type~$A_n$."

\begin{figure}[ht]
\begin{center}
\setlength{\unitlength}{1.5pt}
\begin{picture}(60,66)(0,-5)
\thicklines
  \put(0,20){\line(1,2){20}}
  \put(0,20){\line(1,-1){20}}
  \put(0,20){\line(3,1){60}}
  \put(20,0){\line(0,1){60}}
  \put(20,0){\line(1,1){40}}
  \put(20,60){\line(2,-1){40}}

  \put(20,0){\circle*{1}}
  \put(20,60){\circle*{1}}
  \put(0,20){\circle*{1}}
  \put(60,40){\circle*{1}}

\put(42,17){\makebox(0,0){$a$}}
\put(24,42){\makebox(0,0){$e$}}
\put(34,36){\makebox(0,0){$e'$}}
\put(42,54){\makebox(0,0){$b$}}
\put(7,42){\makebox(0,0){$c$}}
\put(7,7){\makebox(0,0){$d$}}


\end{picture}

\end{center}
\caption{Ptolemy relations}
\label{fig:type-a-exch}
\end{figure}
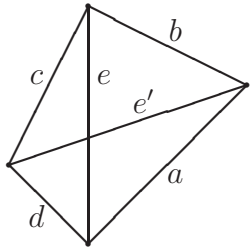

We define the commutative algebra $\AA_n$ explicitly by generators and
relations.
The algebra $\AA_n$ has $\binom{n+3}{2}$ generators~$x_a$, where
the index~$a$ runs over all sides and diagonals of a regular
$(n+3)$-gon~$\mathbf{P}_{n+3}$.
Each defining relation for~$\AA_n$ corresponds to a choice of a
quadruple of vertices in~$\mathbf{P}_{n+3}$.
Let $a$, $b$, $c$, and~$d$ be consecutive sides of the
corresponding convex quadrilateral, and let $e$ and~$e'$
be its diagonals; see Figure~\ref{fig:type-a-exch}.
We then write the relation
\begin{equation}
\label{eq:ptolemy-1}
x_e\,x_{e'}=x_a\,x_c+x_b\,x_d\,.
\end{equation}
Thus, our presentation of $\AA_n$ involves
$\binom{n+3}{4}$ defining relations.

We refer to the relations \eqref{eq:ptolemy-1} as
\emph{Ptolemy relations} because of the classical
Ptolemy Theorem which asserts that in an inscribed
quadrilateral, the sum of the products of the two pairs of
opposite sides equals the product of the two diagonals.
One can pass from Euclidean to hyperbolic
geometry, where an analogue of  \eqref{eq:ptolemy-1}
holds for the exponentiated hyperbolic distances between horocycles drawn
around vertices of a polygon with geodesic sides and cusps at the
vertices (the ``Penner coordinates'' on the corresponding
\emph{decorated Teichm\"uller space}~\cite{penner}).
This observation leads to one of many constructions of cluster algebras
arising in the context of Teichm\"uller theory \cite{fock-goncharov1,
  gsv2, thurston}.

\medskip

The algebra~$\AA_n$ has the following closely related models:
\begin{itemize}
\item
the ring of polynomial $SL_2$-invariants of an
$(n+3)$-tuple of points in~$\CC^2$;
\item
the homogeneous coordinate ring of the Grassmannian
${\rm Gr}_{2,n+3}$ of \hbox{2-dimensional} subspaces in $\CC^{n+3}$
with respect to its Pl\"ucker embedding; 
\item
the ring $\CC[X_n]$, where $X_n$ is the variety of all decomposable
bivectors in~$\bigwedge^2 \CC^{n+3}$.
\end{itemize}
The equivalence of the last two models is trivial.
To see their equivalence to the first model, think of an
$(n+3)$-tuple of points in $\CC^2$ as the columns of a
$2 \times (n+3)$ matrix~$Z$; the row space of $Z$ can then be viewed
as a point of ${\rm Gr}_{2,n+3}$.
The $2 \times 2$ minors $P_{ij}$ of~$Z$ are called
\emph{Pl\"ucker coordinates}; thus, a generic point of~$X_n$
has the form
$$x = \sum_{i < j} P_{ij}\, e_i \wedge e_j,$$
where $e_1, \dots, e_{n+3}$ is the standard basis of~$\CC^{n+3}$.
To identify $\AA_n$ with $\CC[X_n]$, label the vertices of
$\mathbf{P}_{n+3}$ clockwise, and send the generator
$x_a$ corresponding to a chord $a = [i,j]$ to the
Pl\"ucker coordinate $P_{ij}$.
The Ptolemy relations \eqref{eq:ptolemy-1}
then translate into the \emph{Grassmann-Pl\"ucker relations}
\begin{equation}
\label{eq:grassmann-plucker-3term}
P_{ik}\,P_{jl} = P_{ij}\,P_{kl} + P_{il}\,P_{jk},
\end{equation}
for all $1 \leq i < j < k < l \leq n+3$.

One can develop a theory of total positivity in $X_n$ that is
completely parallel to the theory of total positivity in double Bruhat
cells presented in Section~\ref{sec:tp-dbc}.
Let us call a point $x \in X_n$ \emph{totally positive} if all its
Pl\"ucker coordinates are positive real numbers.
We then define \emph{TP-bases} for $X_n$ as $m$-element subsets of
$\AA_n = \CC[X_n]$ satisfying the conditions in
Definition~\ref{def:tp-basis}, with $m = \dim(X_n) = 2n+3$.

In what follows, a \emph{triangulation} of the
polygon~$\mathbf{P}_{n+3}$ always means a triangulation by
non-crossing diagonals.
A nice family of TP-bases for $X_n$ can be constructed as
follows (cf.\ Theorem~\ref{th:tp-basis}).

\begin{theorem}[\cite{ca2}]
\label{th:tp-basis-triangulation}
Every triangulation~$T$ of $\mathbf{P}_{n+3}$ gives rise to a
\linebreak[3]
TP-basis $\tilde \xx(T)$ for~$X_n$ which consists of
the $2n+3$ generators $x_a$ corresponding to the sides and diagonals of~$T$.
\end{theorem}

\begin{remark}
\label{rem:invariant-theory-connection}
The collections $\tilde \xx(T)$ have already appeared in classical
$19^{\rm th}$ century
literature on invariant theory; for invariant-theoretic
connections and applications, see~\cite{kungrota,stur}.
In particular, it is known~\cite{kungrota,stur} that the
monomials in the generators $x_a$ which do not involve
diagonals crossing each other form a linear basis in~$\AA_n$.
We shall later discuss a far-reaching generalization of this
result in the context of cluster algebras.
\end{remark}

\begin{figure}[ht]
\begin{center}
\setlength{\unitlength}{1.6pt}
\begin{picture}(60,60)(0,0)
\thicklines
  \put(0,20){\line(1,2){20}}
  \put(0,20){\line(1,-1){20}}
  \put(20,0){\line(0,1){60}}
  \put(20,0){\line(1,1){40}}
  \put(20,60){\line(2,-1){40}}

  \put(20,0){\circle*{1}}
  \put(20,60){\circle*{1}}
  \put(0,20){\circle*{1}}
  \put(60,40){\circle*{1}}

\put(24,32){\makebox(0,0){$e$}}

\put(7,8){\makebox(0,0){$a$}}
\put(43,18){\makebox(0,0){$b$}}
\put(41,53){\makebox(0,0){$c$}}
\put(6,41){\makebox(0,0){$d$}}
\end{picture}
\begin{picture}(40,66)(0,0)
\put(20,30){\makebox(0,0){$\longrightarrow$}}
\end{picture}
\begin{picture}(60,66)(0,0)
\thicklines
  \put(0,20){\line(1,2){20}}
  \put(0,20){\line(1,-1){20}}
  \put(0,20){\line(3,1){60}}
  \put(20,0){\line(1,1){40}}
  \put(20,60){\line(2,-1){40}}

  \put(20,0){\circle*{1}}
  \put(20,60){\circle*{1}}
  \put(0,20){\circle*{1}}
  \put(60,40){\circle*{1}}

\put(29,34){\makebox(0,0){$e'$}}

\put(7,8){\makebox(0,0){$a$}}
\put(43,18){\makebox(0,0){$b$}}
\put(41,53){\makebox(0,0){$c$}}
\put(6,41){\makebox(0,0){$d$}}

\end{picture}

\end{center}
\caption{A diagonal flip}
\label{fig:flip}
\end{figure}
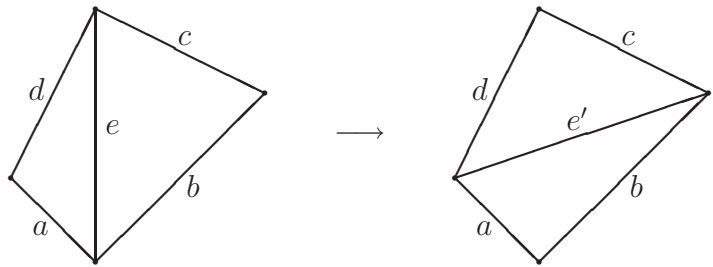

The combinatorics of the family of TP-bases described in
Theorem~\ref{th:tp-basis-triangulation} is encoded by the
\emph{exchange graph} whose vertices are all the triangulations
of~$\mathbf{P}_{n+3}$, and whose edges correspond to the diagonal
\emph{flips}.
Each flip removes a diagonal to create a quadrilateral, then
replaces it with another diagonal of the
same quadrilateral.
See Figure~\ref{fig:flip}.

\pagebreak[3]

\begin{figure}[ht!]
\centerline{\epsfbox{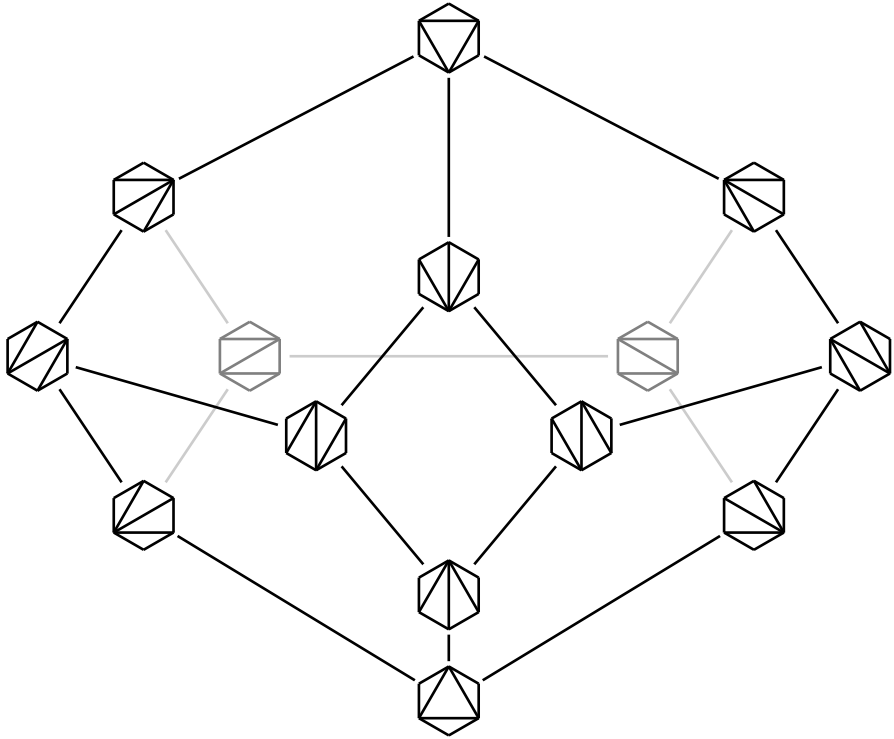}}
\caption{The 3-dimensional associahedron.}
\label{fig:A3assoc_poly}
\end{figure}

Since every triangulation of $\mathbf{P}_{n+3}$ involves~$n$ diagonals,
the exchange graph is $n$-regular, i.e., every vertex has degree~$n$.
The exchange graph is also connected: it is well known that
any two triangulations of~$\mathbf{P}_{n+3}$ are related by a sequence
of flips.
This exchange graph is the $1$-skeleton of an $n$-dimensional
\emph{associahedron}, or Stasheff's polytope~\cite{stasheff}.
See Figure~\ref{fig:A3assoc_poly}.
A generalization of this construction, to be presented in
Section~\ref{sec:gen-asso}, will identify exchange graphs of cluster
algebras of ``finite type'' with the $1$-skeleta of \emph{generalized
  associahedra}, a family of convex polytopes associated with finite
root systems.

For a given triangulation~$T$, each of the~$n$ TP-bases $\tilde \xx(T')$
associated with the neighbors of~$T$ in the exchange graph is obtained
from $\tilde \xx(T)$ by an exchange analogous to \eqref{eq:Fik}:
\begin{equation}
\label{eq:exchange-T-T'}
\tilde \xx(T')  = \tilde \xx(T) - \{x_e\} \cup \{x_{e'}\},
\end{equation}
whenever~$T'$ is obtained from~$T$ by the flip in
Figure~\ref{fig:flip}.
Here the participating generators $x_e$ and $x_{e'}$
satisfy the Ptolemy relation \eqref{eq:ptolemy-1}.

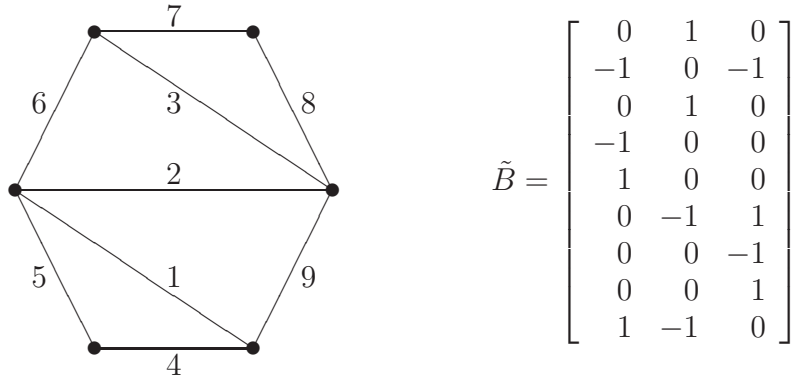
\begin{figure}[ht]
\begin{center}
\setlength{\unitlength}{3pt}
\begin{picture}(100,45)(0,-2)

  \put(0,20){\circle*{1.5}}
  \put(40,20){\circle*{1.5}}
  \put(10,0){\circle*{1.5}}
  \put(10,40){\circle*{1.5}}
  \put(30,0){\circle*{1.5}}
  \put(30,40){\circle*{1.5}}

\put(0,20){\line(1,2){10}}
\put(0,20){\line(1,-2){10}}
\put(40,20){\line(-1,2){10}}
\put(40,20){\line(-1,-2){10}}
\put(10,0){\line(1,0){20}}
\put(10,40){\line(1,0){20}}

\put(0,20){\line(1,0){40}}

\put(0,20){\line(3,-2){30}} \put(10,40){\line(3,-2){30}}

\put(20,9){\makebox(0,0){$1$}}
\put(20,22){\makebox(0,0){$2$}}
\put(20,31){\makebox(0,0){$3$}}

\put(20,-2){\makebox(0,0){$4$}}
\put(20,42){\makebox(0,0){$7$}}
\put(3,9){\makebox(0,0){$5$}}
\put(37,9){\makebox(0,0){$9$}}
\put(3,31){\makebox(0,0){$6$}}
\put(37,31){\makebox(0,0){$8$}}

\put(60,20){{$
\tilde B=
\left[\begin{array}{rrr}
0 & 1 & 0\\
-1& 0 & -1\\
0 & 1 & 0\\
-1& 0 & 0\\
1 & 0 & 0\\
0 &-1 & 1\\
0 & 0 &-1\\
0 & 0 & 1\\
1 &-1 & 0
\end{array}\right]
$}}
\end{picture}
\end{center}
\caption{Matrix $\tilde B$ for a triangulation}
\label{B-tilde-example}
\end{figure}

We now make the following crucial observation:
the~$n$ Ptolemy relations that relate $\tilde \xx(T)$ to its
neighbors can be written in the same form as the exchange relations
\eqref{eq:exchange-f} from Section~\ref{sec:tp-dbc};
that is, they can be encoded by an
$m \times n$ integer matrix $\tilde B = \tilde B(T)$ (as above, $m = 2n+3$).
To be more precise, let $a_1,\dots,a_m$ denote the sides and diagonals
making up the triangulation~$T$.
(We thus fix a labeling of these segments by the numbers $1,\dots,m$.)
Let $\ex = \ex_T \subset [1,m]$ be the $n$-element subset of indices
which label the diagonals of~$T$; as above, we refer to the
elements of $\ex$ as \emph{exchangeable} indices.
Now let $\tilde B = \tilde B(T)=(b_{ij})$ be the $m \times n$ integer matrix
with rows labeled by~$[1,m]$, columns labeled by $\ex$,
and matrix entries given by
\begin{equation}
\label{eq:BT}
b_{ij}=
\begin{cases}
1 & \text{if $a_i$ and $a_j$ are sides in some triangle of~$T$,}\\
  & \text{\ \ with $a_j$ following~$a_i$ in the clockwise order;}\\
-1 & \text{if the same holds, with the counter-clockwise order;}\\
0 & \text{otherwise.}
\end{cases}
\end{equation}
An example is shown in Figure~\ref{B-tilde-example}.

We abbreviate $x_k = x_{a_k}$.
An easy inspection shows that the Ptolemy relation
\eqref{eq:ptolemy-1} between $x_k = x_e$ and $x'_k = x_{e'}$
can indeed be written in the form completely analogous to~\eqref{eq:exchange-f}:
\begin{equation}
\label{eq:exchange-ptolemy}
x'_k = \frac
{\prod\limits_{b_{ik}>0} x_i^{b_{ik}}+
\prod\limits_{b_{ik}<0} x_i^{-b_{ik}}}{x_k} \, .
\end{equation}

We see that the algebraic relationships of the TP-basis $\tilde
\xx(T)$ with its neighbors in the exchange graph are very similar
to the ones involving the TP-basis $F_\ii$ in
Section~\ref{sec:exchanges-double-cells}. There is however an
important difference. In the case of~$\AA_n$, the family of
TP-bases associated with all triangulations of $\mathbf{P}_{n+3}$
is closed under the exchanges, whereas in the case of double
Bruhat cells, this is not the case (see Example~\ref{ex:sl3-exchanges}).
This raises the following question, whose answer
constituted the decisive step towards introducing cluster algebras:
is there a mechanism for extending the exchanges from
an ``initial" TP-basis $F_\ii$ to all its $n$ neighbors~$F_{\ii;k}$,
then to the $n$ neighbors of~$F_{\ii;k}$, etc.?
More concretely, is there a natural general rule that determines the matrix
$\tilde B'$ associated with a neighboring TP-basis 
from the matrix $\tilde B$ associated with a given TP-basis?
Such a rule would enable us to recursively ``propagate'' away from the
initial TP-basis, each time constructing a new basis using the
exchange relations involving the entries of the current matrix~$\tilde
B$, then re-calculating $\tilde B$ according to the stated rule.

The case of triangulations provides a crucial hint for the
solution: it turns out that the matrix $\tilde B(T')$
can be recovered from $\tilde B(T)$ by a purely algebraic
procedure which we call \emph{matrix mutation}\footnote{Our historic
  account becomes a bit revisionist
at this point: in fact, we discovered matrix mutations
\emph{before} considering the algebras~$\AA_n\,$.}.

\begin{definition}[{\cite[Definition~4.2]{ca1}}]
\label{def:matrix-mutation}
We say that an
$m \times n$ integer matrix $\tilde B'$ is obtained from $\tilde B$ by
\emph{matrix mutation} in direction~$k \in \ex$,
and write $\tilde B' = \mu_k (\tilde B)$
if the entries of $\tilde B'$ are given by
\begin{equation}
\label{eq:matrix-mutation}
b'_{ij} =
\begin{cases}
-b_{ij} & \text{if $i=k$ or $j=k$;} \\[.05in]
b_{ij} + \displaystyle\frac{|b_{ik}| b_{kj} +
b_{ik} |b_{kj}|}{2} & \text{otherwise.}
\end{cases}
\end{equation}
\end{definition}

The following proposition can be verified by direct inspection.

\begin{proposition}
\label{pr:mut-tri}
Suppose that a triangulation $T'$ is obtained from $T$ by
a flip replacing a diagonal~$a_k$ by a diagonal~$a'_k$.
(The labeling of all other diagonals and sides remains
unchanged.)
Then $\tilde B(T')= \mu_k(\tilde B(T))$.
\end{proposition}

\section{Cluster algebra fundamentals 
}
\label{sec:fundam}

Motivated by the above considerations, we are now ready to develop
the axiomatic foundations of the cluster algebra theory.

\subsection{Upper cluster algebras}

We start by introducing
\emph{upper cluster algebras} (of geometric type),
which first appeared in~\cite{ca3}.

Let~$m$ and~$n$ be two positive integers such that~$m \geq n$.
Let~$\FFcal$ be the field of rational
functions over~$\QQ$ in~$m$ independent (commuting) variables.
Once and for all, we fix an $n$-element subset $\ex \subset [1,m]$.

\begin{definition}
\label{def:seed}
A \emph{seed} (of geometric type) in $\FFcal$ is a pair
$(\tilde \xx, \tilde B)$, where
\begin{itemize}

\item $\tilde \xx = \{x_1, \dots, x_m\}$ is an
algebraically independent subset of~$\FFcal$ which generates~$\FFcal$;

\item $\tilde B=(b_{ij})$ is an $m \times n$ integer matrix with rows
labeled by $[1,m]$
and columns labeled by~$\ex$, such that
\begin{itemize}
\item[$\circ$]
the $n \times n$ submatrix $B\!=\!(b_{i,j})_{i,j\in\ex}$ is
\emph{skew-symmetrizable}:
$d_i b_{ik} = -d_k b_{ki}$ for some positive
integers~$d_i \,\, (i, k \in \ex)$;
\item[$\circ$]
$\tilde B$ has full rank~$n$.
\end{itemize}
\end{itemize}
The seeds are defined up to a relabeling of elements of $\tilde \xx$
together with the corresponding relabeling of rows and columns of~$\tilde B$.
\end{definition}

\begin{example}
In the notation of Section~\ref{sec:tp-dbc},
let $\ii$ be a reduced word for $(u,v)\in W\times W$.
As shown in \cite[Proposition~2.6]{ca3},
the matrix $\tilde B =\tilde B(\ii)$ given by \eqref{eq:tildeB-entries}
satisfies the conditions in Definition~\ref{def:seed}.
Thus, $(F_\ii,\tilde B(\ii))$ is a seed inside
the field of rational functions $\FFcal=\QQ(G^{u,v})$.
\end{example}

\begin{example}
In the notation of Section~\ref{sec:triangulations},
let $T$ be a triangulation of the regular $(n+3)$-gon
$\mathbf{P}_{n+3}$.
One can show that the matrix $\tilde B =
\tilde B(T)$ defined by \eqref{eq:BT}
satisfies the conditions in Definition~\ref{def:seed}.
Thus, $(\tilde \xx(T),\tilde B(T))$ is a seed inside~$\FFcal$,
the field  of fractions of the algebra~$\AA_n$.
\end{example}

Let $(\tilde \xx, \tilde B)$ be a seed in a field~$\FFcal$.
By analogy with \eqref{eq:exchange-f} and
\eqref{eq:exchange-ptolemy},
for every~$k \in \ex$, we set
\begin{equation}
\label{eq:exchange}
x'_k = \frac
{\prod\limits_{b_{ik}>0} x_i^{b_{ik}}+
\prod\limits_{b_{ik}<0} x_i^{-b_{ik}}}{x_k} \in \FFcal
\end{equation}
and
\begin{equation}
\label{eq:adjacent-cluster}
\tilde \xx_k = \tilde \xx - \{x_k\} \cup \{x'_k\} \, .
\end{equation}
It is easy to check that $\tilde \xx_k$ is again a set of algebraically
independent generators for~$\FFcal$.

\begin{definition}
\label{def:upper-cluster-algebra}
The \emph{upper cluster algebra}
$\Upper(\tilde \xx, \tilde B) = \Upper(\tilde B)$
is the subring of the ambient field~$\FFcal$ given by
\begin{equation}
\label{eq:upper}
\Upper(\tilde \xx, \tilde B) = \Upper(\tilde B)
 = \ZZ[\tilde \xx^{\pm 1}] \cap \bigcap_{k \in \ex}
\ZZ[\tilde \xx_k^{\pm 1}],
\end{equation}
where $\ZZ[\tilde \xx^{\pm 1}]$ (resp., $\ZZ[\tilde \xx^{\pm 1}_k]$)
stands for the ring of
Laurent polynomials with integer coefficients in the variables from~$\tilde \xx$
(resp.,~$\tilde \xx_k$).
\end{definition}

In this terminology, Corollary~\ref{cor:CGuv} takes the following
form (cf.~\cite[Theorem~2.10]{ca3}):

\begin{corollary}
\label{cor:CGuv-upper-cluster}
The coordinate ring $\CC[G^{u,v}]$ of any double Bruhat cell
is naturally isomorphic to the complexification of the upper
cluster algebra $\Upper(\tilde B(\ii))$
associated with the matrix $\tilde B(\ii)$
determined by an arbitrary reduced word $\ii$  for $(u,v)$.
\end{corollary}

To make the abstract notion of an upper cluster algebra useful,
we need to develop a better understanding of the structure of these algebras.
This does not seem to be an easy task.
In particular, it is not at all clear how big the intersection
of Laurent polynomial rings in \eqref{eq:upper} can be.
The following problem raised in \cite[Problem~1.27]{ca3}
is still open:

\begin{problem}
\label{prob:finite-gen-upper}
For which matrices~$\tilde B$ is the upper cluster algebra~$\Upper(\tilde B)$
finitely generated?
\end{problem}

\subsection{Seed mutations and cluster algebras}

Prompted by Proposition~\ref{pr:mut-tri}, we will now introduce
the machinery of seed mutations,
which will provide the means to generate a lot of---sometimes
infinitely many---new elements
in~$\Upper(\tilde x,\tilde B)$ 
from the initial seed~$(\tilde x,\tilde B)$.
This will in particular lead to partial answers to
Problem~\ref{prob:finite-gen-upper}
(see Theorem~\ref{th:acyclic}).

We retain the terminology and notation in
Definition~\ref{def:seed}.
In addition, we denote
$\xx = \{x_j: j \in \ex\} \subset \tilde \xx$ and
$\cc = \tilde \xx - \xx$
(cf.~\eqref{eq:cc-bdc}).
We refer to the elements of $\ex$ as \emph{exchangeable} indices,
and to $\xx$ as a \emph{cluster}.
The $n \times n$ submatrix $B$ of $\tilde B$ with rows and columns
labeled by~$\ex$ is called the \emph{principal part} of~$\tilde B$.
We also call $B$ the \emph{exchange matrix} of the seed
$(\tilde \xx, \tilde B)$.

Recall the notion of matrix mutation introduced in
Definition~\ref{def:matrix-mutation}.
We note that if the principal part of $\tilde B$ is~$B$,
then the principal part of $\mu_k(\tilde B)$ is~$\mu_k(B)$.
One can show (see, e.g., \cite[Proposition~2.3]{qca})
that matrix mutations preserve the restrictions imposed on $\tilde B$
in Definition~\ref{def:seed}.
Consequently, the following notion is well defined.

\begin{definition}
\label{def:seed-mutation}
Let $(\tilde \xx, \tilde B)$ be a seed in $\FFcal$.
For any exchangeable index $k$, the \emph{seed mutation} in direction~$k$
transforms $(\tilde \xx, \tilde B)$ into the seed
$\mu_k(\tilde \xx, \tilde B)=(\tilde \xx', \tilde B')$, where
$\tilde \xx' = \tilde \xx_k$ is given by
\eqref{eq:exchange}--\eqref{eq:adjacent-cluster}, and $\tilde B' =
\mu_k(\tilde B)$.
\end{definition}

It is trivial to check that matrix mutations are involutive,
as are seed mutations: in the situation of
Definition~\ref{def:seed-mutation} we have
$\mu_k(\tilde B') = \tilde B$ and
$\mu_k (\tilde \xx', \tilde B') = (\tilde \xx, \tilde B)$.
Therefore, we can define an equivalence relation on seeds as follows:
we say that two seeds $(\tilde \xx, \tilde B)$ and $(\tilde \xx', \tilde B')$
are mutation-equivalent, and write
$(\tilde \xx, \tilde B) \sim (\tilde \xx', \tilde B')$, if $(\tilde \xx', \tilde B')$
can be obtained from $(\tilde \xx, \tilde B)$ by a sequence of seed mutations.
All seeds $(\tilde \xx', \tilde B')$ mutation-equivalent to
a given seed $(\tilde \xx, \tilde B)$ share the same set
$\cc = \tilde \xx - \xx$.

\begin{remark}
Recall that our motivating example of an upper cluster algebra is
the coordinate ring $\QQ[G^{u,v}]$ of a double Bruhat cell,
see Corollary~\ref{cor:CGuv-upper-cluster}.
The mutation procedure makes it obvious that, for every seed
$(\tilde \xx, \tilde B)$ mutation-equivalent to
$(F_\ii, \tilde B(\ii))$, the subset
$\tilde \xx \subset \QQ[G^{u,v}]$ satisfies conditions (i) and (iii)
in Definition~\ref{def:tp-basis}, and therefore provides global
coordinates in the totally positive variety $G^{u,v}_{> 0}$.
\end{remark}

Somewhat surprisingly, seed mutations leave the associated upper cluster
algebra invariant:

\begin{theorem}[{\cite[Theorem~1.5]{ca3}}]
\label{th:upper-mutation-invariant}
If $(\tilde \xx,\tilde B)\sim (\tilde \xx',\tilde B')$,
then $\Upper(\tilde \xx,\tilde B)=\Upper(\tilde \xx',\tilde B')$.
\end{theorem}

In view of Theorem~\ref{th:upper-mutation-invariant}, we can use the notation
$\Upper(\mathcal{S}) = \Upper(\tilde \xx,\tilde B)$, where
$\mathcal{S}$ is the mutation-equivalence class of a seed
$(\tilde \xx,\tilde B)$.


Now everything is in place for defining cluster algebras.
For a mutation equivalence class of seeds~$\mathcal{S}$,
we denote by $\mathcal{X} = \mathcal{X}(\mathcal{S})$
the union of clusters of all seeds in~$\mathcal{S}$.
In other words, $\mathcal{X}$ is the union of all elements of~$\FFcal$
that appear in all seeds of~$\mathcal{S}$, except that we do not
include the elements of~$\cc$.
We refer to the elements of $\mathcal{X}$ as \emph{cluster variables}.

\begin{definition}
\label{def:cluster-algebra}
Let $\mathcal{S}$ be a mutation equivalence class of seeds.
The \emph{cluster algebra} $\AA(\mathcal{S})$ associated with $\mathcal{S}$
is the $\ZZ[\cc^{\pm 1}]$-subalgebra\footnote{
An alternative definition uses $\ZZ[\cc]$ instead of $\ZZ[\cc^{\pm
  1}]$
as a ground ring for~$\AA(\mathcal{S})$; cf.\
Section~\ref{sec:triangulations}.
More generally, one can use any ground ring sandwiched between
$\ZZ[\cc]$ and $\ZZ[\cc^{\pm1}]$.
}
of the ambient field $\FFcal$
generated by all cluster variables:
\[
\AA(\mathcal{S}) \stackrel{\rm def}{=} \ZZ[\cc^{\pm 1}, \mathcal{X}].
\]
The cardinality~$n$ of every cluster is called the
\emph{rank} of $\AA(\mathcal{S})$.
\end{definition}

Since $\mathcal{S}$ is uniquely determined by an arbitrary seed
$(\tilde \xx, \tilde B)\in\mathcal{S}$, we sometimes
denote $\AA(\mathcal{S})$ as $\AA(\tilde \xx,\tilde B)$,
or simply as~$\AA(\tilde B)$,
because $\tilde B$ determines this algebra uniquely up to an
automorphism of~$\FFcal$.

The following result first appeared in~\cite[Theorem~3.1]{ca1}.
A new proof based on Theorem~\ref{th:upper-mutation-invariant} was
given in~\cite{ca3}.

\begin{theorem}[Laurent phenomenon]
\label{th:A-in-U}
A cluster algebra~$\AA(\mathcal{S})$ is contained in
the corresponding upper cluster algebra~$\Upper(\mathcal{S})$.
Equivalently,
$\AA(\mathcal{S}) \ \subset \ \ZZ[\tilde \xx^{\pm 1}]$ for every seed
$(\tilde \xx, \tilde B)$ in~$\mathcal{S}$.
That is, every element of $\AA(\mathcal{S})$ is an integer
Laurent polynomial in the variables from~$\tilde \xx$.
\end{theorem}

It is worth mentioning that Theorems~\ref{th:upper-mutation-invariant}
and~\ref{th:A-in-U}
carry over to the quantum setting developed in~\cite{qca}.

We next present some open problems and conjectures on cluster algebras.
The following question was raised in \cite[Problem~1.25]{ca3}.

\begin{problem}
\label{problem:when-upper=cluster}
When does a cluster algebra~$\AA(\mathcal{S})$ coincide with the
corresponding upper cluster algebra~$\Upper(\mathcal{S})$?
\end{problem}

\pagebreak[3]

A sufficient condition for the equality
$\AA(\mathcal{S}) = \Upper(\mathcal{S})$ was found in~\cite{ca3}.
Following~\cite{ca3}, let us call a seed~$(\tilde \xx,\tilde B)$
\emph{acyclic} if there is a linear ordering of~$\ex$ such that
$b_{ij} \geq 0$ for all $i,j \in \ex$ with $i < j$.
(Since the principal part of~$\tilde B$ is skew-symmetrizable,
we also have $b_{ji} \leq 0$.)

\begin{theorem}[{\cite[Theorem~1.18, Corollary~1.19]{ca3}}]
\label{th:acyclic}
If~$\mathcal{S}$ contains  an acyclic seed,
then $\AA(\mathcal{S}) = \Upper(\mathcal{S})$.
Moreover, if $(\tilde \xx,\tilde B)$ is an acyclic seed, then
\[
\AA(\tilde \xx,\tilde B) = \Upper(\tilde \xx,\tilde B)
= \ZZ[\cc^{\pm 1}, x_k, x'_k \, (k \in \ex)],
\]
where the $x'_k$ is given by~\eqref{eq:exchange}.
In particular, this algebra is finitely generated.
\end{theorem}

An important combinatorial invariant of a cluster algebra
$\AA(\mathcal{S})$ is its \emph{exchange graph}.
The vertices of this graph correspond to the seeds in~$\mathcal{S}$,
and the edges correspond to seed mutations.
By definition, the exchange graph is connected and $n$-regular;
that is, every vertex has degree~$n$.
Besides these properties, not much is known about exchange
graphs in general.
Here are a few conjectures (some of them already
appeared in~\cite{ca2}).

\begin{conjecture}
\label{con:exchange-general}
{\ }
\begin{enumerate}
\item
The exchange graph of an algebra $\AA(\tilde \xx,\tilde B)$
depends only on the principal part~$B$ of the matrix~$\tilde B$.
\item
Every seed in~$\mathcal{S}$ is uniquely determined by its cluster;
thus, the vertices of the exchange graph can be identified with the
clusters,
with two clusters adjacent if and only if their
intersection has cardinality~$n-1$.
\item
For any cluster variable $x$, the seeds
from~$\mathcal{S}$ whose clusters contain~$x$ form a
{connected} subgraph of the exchange graph.
\item
The seeds from~$\mathcal{S}$ whose exchange matrix~$B$ is acyclic
form a {connected} subgraph (possibly empty) of the exchange graph.
\end{enumerate}
\end{conjecture}


\subsection{Cluster monomials and positive cones}

\begin{definition}
\label{def:cluster-monomial}
A \emph{cluster monomial} in a cluster algebra $\AA(\mathcal{S})$
is a monomial in cluster variables all of which
belong to the same cluster.
\end{definition}

Inspired by Remark~\ref{rem:invariant-theory-connection}, we
propose the following conjecture.

\begin{conjecture}
\label{con:cluster-monomials-lin-indep}
The cluster monomials are
linearly independent over~$\ZZ[\cc^{\pm 1}]$.
\end{conjecture}

As labels for cluster monomials, we would like to use the
denominators in the Laurent expansion with respect to a given
cluster. To be more specific, let us fix a cluster~$\xx$ from some
seed in~$\mathcal{S}$, and choose a numbering of its elements:
$\xx = \{x_1, \dots, x_n\}$. Every nonzero element $y \in
\Upper(\mathcal{S})$ can be uniquely written as
\begin{equation}
\label{eq:Laurent-normal-form}
y = \frac{P(x_1, \dots, x_n)}{x_1^{d_1} \cdots x_n^{d_n}} \, ,
\end{equation}
where $P(x_1, \dots, x_n)$ is a polynomial with coefficients in~$\ZZ[\cc^{\pm 1}]$
which is not divisible by any cluster variable~$x_i\in\xx$.
We denote
\begin{equation}
\label{eq:denominator-vector}
\delta(y) = \delta_\xx(y) = (d_1, \dots, d_n) \in \ZZ^n,
\end{equation}
and call the integer vector~$\delta(y)$ the \emph{denominator vector}
of~$y$ with respect to the cluster~$\xx$.
(One could also call $\delta(y)$ the \emph{tropicalization} of $y$
with respect to~$\xx$.)
For instance, the elements of~$\xx$ and their immediate exchange
partners have denominator vectors
\begin{equation}
\label{eq:denom-xj}
\delta (x_j) = - e_j, \quad \delta (x'_j) = e_j \quad (j \in [1,n]),
\end{equation}
where $e_1, \dots, e_n$ are the standard basis vectors in~$\ZZ^n$.
Note also that the map $y \mapsto \delta(y)$ has the
following valuation property:
\begin{equation}
\label{eq:delta-multiplicative}
\delta(yz) = \delta(y) + \delta(z) \,.
\end{equation}

\begin{conjecture}
\label{con:denominators}
Different cluster monomials have different denominator vectors with
respect to a given cluster.
\end{conjecture}

Our last group of questions/conjectures concerns \emph{positivity} in
cluster algebras.
We define the \emph{positive cones} $\Upper(\mathcal{S})_{\geq 0}$
and $\AA(\mathcal{S})_{\geq 0}$ by setting
\begin{align*}
\Upper(\mathcal{S})_{\geq 0}
&= \bigcap_{(\tilde \xx,\tilde B) \in \mathcal{S}}
\ZZ_{\geq 0}[\tilde \xx^{\pm 1}], \\[.1in]
\AA(\mathcal{S})_{\geq 0}
&= \AA(\mathcal{S}) \ \cap \ \Upper(\mathcal{S})_{\geq 0}\, .
\end{align*}
That is, $\Upper(\mathcal{S})_{\geq 0}$ (resp., $\AA(\mathcal{S})_{\geq 0}$)
consists of all elements in $\Upper(\mathcal{S})$
(resp., in~$\AA(\mathcal{S})$) whose Laurent expansion
in terms of the variables from any seed in $\mathcal{S}$ has nonnegative
coefficients.
A nonzero element of a positive cone will be referred to as \emph{positive}.

It would be very interesting to find concrete
descriptions of the positive cones $\Upper(\mathcal{S})_{\geq 0}$
and $\AA(\mathcal{S})_{\geq 0}$.
To be more specific: let us call a positive element (in either of the
positive cones) \emph{indecomposable} if it cannot be written as
a sum of two positive elements.

\begin{problem}
\label{prob:indec-positive}
Describe the set of all indecomposable elements
in either of the cones $\Upper(\mathcal{S})_{\geq 0}$
and $\AA(\mathcal{S})_{\geq 0}$.
\end{problem}

\pagebreak[2]

\begin{conjecture}
\label{con:cluster-positive}
{\ }
\begin{enumerate}
\item
\cite{ca1}
Every cluster variable $x\in\mathcal{X}$ is a positive element of the
cluster algebra~$\AA(\mathcal{S})$.
That is, for every seed $(\tilde \xx,\tilde B)$ in~$\mathcal{S}$,
the Laurent expansion of $x$ in terms of  the variables from~$\tilde
\xx$ has nonnegative coefficients.
\item
Every cluster monomial is an indecomposable positive element.
\end{enumerate}
\end{conjecture}

\subsection{Example: cluster algebras of rank 2}
Let us illustrate the above definitions, results and conjectures
by the following example.
Let $m = n = 2$, and let the ambient field
$\mathcal{F}$ be the field of rational functions
$\QQ(y_1, y_2)$ in two independent variables $y_1$ and $y_2$.
For every pair of positive integers~$(b,c)$, let $\AA(b,c)$ be the
cluster algebra associated with the initial seed
$(\{y_1, y_2\}, B)$, where $B = \tilde B$ is the $2 \times 2$
matrix given by
$$B \!= \! \tilde B \! = \! \bmat{0}{b}{-c}{0} \, .$$
The definitions readily imply that the cluster variables in
$\AA(b,c)$ are the elements $y_t \in \FFcal$, for $t \in \ZZ$, defined
recursively by
\[
y_{t-1} y_{t+1} =
\left\{\!
\begin{array}{ll}
y_t^b + 1 & {\rm if~} t {\rm~is \, odd}; \\[.2in]
y_t^c + 1 & {\rm if~} t {\rm~is \, even}.
\end{array}
\right.
\]
Thus, $\AA(b,c)$ is simply the subring of $\FFcal=\QQ(y_1, y_2)$ generated by
all the elements~$y_t$.
The clusters are of the form $\{y_t, y_{t+1}\}$, $t \in \ZZ$.
The Laurent phenomenon asserts that
$y_t \in \ZZ[y_1^{\pm 1}, y_2^{\pm 1}]$ for all $t \in \ZZ$.
In view of Theorem~\ref{th:acyclic}, we have a sharper statement:
$$\AA(b,c) = \bigcap_{t=0}^2 \ZZ[y_t^{\pm 1}, y_{t+1}^{\pm 1}]
= \bigcap_{t \in \ZZ} \ZZ[y_t^{\pm 1}, y_{t+1}^{\pm 1}]
= \ZZ[y_0, y_1, y_2, y_3].$$

\begin{theorem}
\label{th:fintype-rank2}
{\rm \cite{ca1}}
The sequence $(y_t)_{t\in\ZZ}$ of the cluster variables in $\AA(b,c)$
is periodic if and only if $bc \leq 3$.
For $bc = 1$ (resp.,~$2$, $3$),
the sequence $(y_t)_{t\in\ZZ}$ has period~$5$ (resp.,~$6$,~$8$).
\end{theorem}

Conjectures~\ref{con:exchange-general}, \ref{con:cluster-monomials-lin-indep}
 and \ref{con:denominators}
are known to hold for any $\AA(b,c)$.
(The first one follows from
the results in \cite{ca1}, while the last two were proved in~\cite{sz}.)
Conjecture~\ref{con:cluster-positive} is still open even in
rank~$2$; it was proved in~\cite{sz} for the special case
$bc \leq 4$.
Furthermore, it was shown in \cite{sz} that for $bc \leq 4$,
the indecomposable positive elements form a $\ZZ$-basis of the cluster algebra.
This raises the following general question.

\begin{problem}
\label{prob:indec-positive-basis}
For which matrices~$\tilde B$
do the indecomposable positive elements form a $\ZZ$-basis
of the corresponding cluster algebra?
\end{problem}

\subsection{Cluster algebras of finite type}
All the problems and conjectures stated above become much more tractable
for the cluster algebras of finite type, defined as follows.

\begin{definition}
\label{def:finite-type}
A cluster algebra $\AA(\mathcal{S})$ is of
\emph{finite type} if the mutation equivalence class~$\mathcal{S}$
consists of finitely many seeds.
\end{definition}

To rephrase, a cluster algebra  is of finite type
if its exchange graph is finite.
In particular, Theorem~\ref{th:fintype-rank2} shows that
a rank~$2$ cluster algebra $\AA(b,c)$ is of finite type if and
only if~$bc \leq 3$.

\begin{theorem}[\cite{ca2}]
\label{th:finite-acyclic}
Every cluster algebra of finite type contains  an acyclic seed, and
therefore satisfies the conclusions in Theorem~\ref{th:acyclic}.
\end{theorem}

\begin{theorem}[\cite{ca2}]
\label{th:exchange-general-finite}
Parts~(1)-(3) of Conjecture~\ref{con:exchange-general} hold
for any cluster algebra~$\AA$ of finite type.
In particular, each seed in~$\AA$ 
is uniquely determined by its cluster.
\end{theorem}

\begin{theorem}[\cite{ca2}]
\label{th:cluster-monomial-indep-finite}
Conjecture~\ref{con:cluster-monomials-lin-indep} holds
for any cluster algebra~$\AA$ of finite type.
\end{theorem}

The proofs of Theorems~\ref{th:finite-acyclic}--
\ref{th:cluster-monomial-indep-finite} are based on the
classification of cluster algebras of finite type given in
\cite{ca2}. Remarkably, this classification (to be discussed in
Section~\ref{sec:finite-type}) turns out to be identical to the
famous Cartan-Killing classification of semisimple Lie algebras
and finite root systems. In the cases where a cluster algebra is
of classical Cartan-Killing type, stronger results can be
obtained.

\begin{theorem}[\cite{ca4}]
\label{th:conjectures-classical}
Conjectures~\ref{con:exchange-general},
\ref{con:cluster-monomials-lin-indep},
\ref{con:denominators} and \ref{con:cluster-positive} hold
for any cluster algebra of classical type.
\end{theorem}

In fact, Theorem~\ref{th:conjectures-classical} can be sharpened
as follows.

\begin{theorem}[\cite{ca4}]
\label{th:conjectures-classical-stronger}
For any cluster algebra~$\AA$ of classical type,
the cluster monomials form a $\ZZ[\cc^{\pm 1}]$-linear basis
of~$\AA$ and a $\ZZ_{\geq 0}[\cc^{\pm 1}]$-linear basis
of the positive cone $\AA_{\geq 0}$.
\end{theorem}

\begin{conjecture}
\label{con:characterizing-finite-type}
Each of the two properties in
Theorem~\ref{th:conjectures-classical-stronger}
characterizes cluster algebras of finite type; that is, it holds
for a cluster algebra~$\AA$ if and only if $\AA$ is of finite type.
\end{conjecture}


\section{Finite type classification and cluster combinatorics}
\label{sec:finite-type}

\subsection{Finite type classification}
We now present a classification of cluster algebras of finite type
given in \cite{ca2}.
First, a few preliminaries.

\begin{figure}[ht]
\begin{center}
\setlength{\unitlength}{2pt}
\begin{picture}(140,34)(-10,-23)
\put(0,0){\line(1,0){120}}
\put(40,0){\line(0,-1){20}}
\put(40,-20){\circle*{2}}
\multiput(0,0)(20,0){7}{\circle*{2}}
\put(0,4){\makebox(0,0){$1$}}
\put(40,-25){\makebox(0,0){$-1$}}
\put(20,4){\makebox(0,0){$-1$}}
\put(40,4){\makebox(0,0){$1$}}
\put(60,4){\makebox(0,0){$-1$}}
\put(80,4){\makebox(0,0){$1$}}
\put(100,4){\makebox(0,0){$-1$}}
\put(120,4){\makebox(0,0){$1$}}
\end{picture}
\end{center}
\caption{Sign function on a Dynkin diagram}
\label{fig:signs-dynkin}
\end{figure}
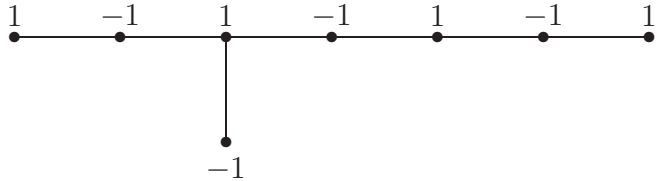

Let $A=(a_{ij})$ be an $n\times n$ Cartan matrix of finite type.
The associated Dynkin diagram is a tree, hence a bipartite graph.
Thus, there is a sign function $\varepsilon: [1,n] \to \{1,-1\}$
such that $a_{ij}< 0 \ \Longrightarrow\
\varepsilon(i)=-\varepsilon(j)$. \linebreak[3]
See Figure~\ref{fig:signs-dynkin}.

We define a skew-symmetrizable  matrix $B(A)\!=\!(b_{ij})$ by
\begin{equation}
\label{eq:B(A)}
b_{ij}=
\begin{cases}
\ \ 0 & \text{if $i=j$;} \\
\varepsilon(i)\, a_{ij} & \text{if $i\neq j$.}
\end{cases}
\end{equation}

\begin{theorem}[\cite{ca2}]
\label{th:fin-type-class}
A cluster algebra $\AA$ is of finite type
if and only if the exchange matrix at some
seed of~$\AA$ is of the form~$B(A)$, where $A$ is a Cartan matrix
of finite type.
Furthermore, the type of~$A$ in the Cartan-Killing
nomenclature is uniquely determined by the cluster
algebra~$\AA$. (It is called the ``cluster type'' of~$\AA$.)
\end{theorem}


\begin{example}
In rank~$2$, the cluster algebras of finite type were classified in
Theorem~\ref{th:fintype-rank2}.
The values of $bc$ equal to 
$1$, $2$, and $3$ correspond to the
cluster types 
$A_2$, $B_2$, and~$G_2$, respectively.
\end{example}

\begin{example}
\label{ex:type-An}
The cluster algebra $\AA_n$ discussed in
Section~\ref{sec:triangulations} has cluster type~$A_n$.
To see this, consider a ``snake" triangulation~$T$ of the regular
$(n+3)$-gon $\mathbf{P}_{n+3}$ (Figure~\ref{B-tilde-example}
illustrates the case $n=3$)
and note that the principal part of the corresponding matrix
$\tilde B(T)$ given by \eqref{eq:BT} is of the form
\eqref{eq:B(A)}, where~$A$ is the Cartan matrix of type~$A_n$.
\end{example}

\begin{remark}
It is a non-trivial problem to decide, given a matrix~$B$, whether the
corresponding cluster algebra is of finite type.
This problem was solved by A.~Seven in~\cite{seven-2}.
An alternative 
criterion has been recently given
in~\cite{bgz}.
\end{remark}

\begin{remark}
\label{rem:dbc-fintype}
It would be interesting to classify the double Bruhat cells such
that the corresponding cluster algebra is of finite type.
Some examples of this kind were given in \cite{ca3}.
For instance, as shown in \cite[Example~2.18]{ca3}, the cluster
algebra structure in $\QQ[SL_3^{w_\circ,w_\circ}]$
(see Examples~\ref{example:SL3-w0w0} and \ref{ex:sl3-exchanges} above)
is of cluster type~$D_4$.
More details are given in Example~\ref{ex:sl3-D4} below.
\end{remark}

Over the last several years, cluster algebra structures have been
uncovered in (homogeneous) coordinate rings of various classical algebraic
varieties.
Some of these cluster algebras have finite type; see
Figure~\ref{fig:cluster-types}.

\begin{figure}[ht]
\[
\begin{array}{ll}
\QQ[\mathrm{Gr}_{2,n+3}]\hspace{1in} & A_n\\[.2in]
\QQ[\mathrm{Gr}_{3,6}] & D_4\\[.1in]
\QQ[\mathrm{Gr}_{3,7}] & E_6\\[.1in]
\QQ[\mathrm{Gr}_{3,8}] & E_8\\[.2in]
\QQ[SL_3/N] & A_1\\[.1in]
\QQ[SL_4/N] & A_3\\[.1in]
\QQ[SL_5/N] & D_6\\[.2in]
\QQ[Sp_4/N] & B_2\\[.2in]
\QQ[SL_2] & A_1\\[.1in]
\QQ[SL_3] & D_4
\end{array}
\]
\caption{Cluster types of some coordinate rings}
\label{fig:cluster-types}
\end{figure}

\pagebreak[3]

In most cases, the symmetry exhibited by the cluster type
of a cluster algebra is not apparent at all from its geometric
realization.
An instance of this phenomenon was demonstrated by
J.~Scott~\cite{scott}, who introduced a natural cluster algebra
structure in the homogeneous coordinate ring of the
Grassmannian~${\rm Gr}_{k,n}$ of $k$-dimensional subspaces
of~$\CC^n$ with respect to its Pl\"ucker embedding.
(Cf.\ also~\cite{gsv1}.)
Among these cluster algebras,
all which are of finite type are listed in Figure~\ref{fig:cluster-types}.

\subsection{Cluster combinatorics and generalized associahedra}
\label{sec:gen-asso}

Let~$\AA$ be a cluster algebra of finite type.
In view of Theorem~\ref{th:exchange-general-finite},
the combinatorics of exchanges and mutations in~$\AA$
is encoded by the \emph{cluster complex} $\Delta(\AA)$,  
a simplicial complex (indeed, a pseudomanifold)
on the set of all cluster variables
whose maximal simplices are the clusters.

For example, the simplices of the cluster complex $\Delta(\AA_n)$
associated with the cluster algebra $\AA_n$ of
Section~\ref{sec:triangulations}
are naturally identified with collections of non-crossing
diagonals of the regular $(n+3)$-gon~$\mathbf{P}_{n+3}$.
In what follows, we extend this description to all cluster
algebras of finite type.

According to Theorem~\ref{th:fin-type-class}, $\AA$ has a seed
$(\tilde \xx, \tilde B)$ with the exchange matrix
$B = B(A)$ given by \eqref{eq:B(A)}, where $A=(a_{ij})$
is a Cartan matrix of finite type (see
\eqref{eq:cartan-transition}).
We refer to such a seed 
as a \emph{distinguished} seed of~$\AA$. 

\begin{theorem}[\cite{ca2}]
\label{th:conjectures-finite}
Conjectures~\ref{con:denominators}
and \ref{con:cluster-positive}(1) hold for a distinguished
seed in~$\AA$.
\end{theorem}

The proof of Conjecture~\ref{con:denominators} (for a distinguished seed)
given in \cite[Theorem~1.9]{ca2} establishes a direct connection between
a cluster algebra of finite type and the corresponding root system.
Let $\Phi$ be a root system with the Cartan matrix~$A$,
and $Q$ the root lattice generated by~$\Phi$.
We identify~$Q$ with $\ZZ^n$ using the
basis~$\Pi=\{\alpha_1,\dots,\alpha_n\}$ of simple roots in~$\Phi$.
Let~$\Phi_{> 0}$ be the set of positive roots associated to~$\Pi$.
In this notation, \cite[Theorem~1.9]{ca2} combined with \cite[Theorem~1.8]{yga}
can be stated as follows.

\begin{theorem}
\label{th:denominator-roots}
The denominator vector parametrization
(see \eqref{eq:Laurent-normal-form}--\eqref{eq:denominator-vector})
with respect to a distinguished cluster~$\xx$ provides a bijection
between the set of cluster variables $\mathcal{X}$ and the set
\[
\Phi_{\geq -1} = \Phi_{> 0} \cup (- \Pi)
\]
of ``almost positive roots."
This parametrization also gives a bijection between the set of all cluster
monomials and the root lattice~$Q = \ZZ^n$.
\end{theorem}

We next restate Theorem~\ref{th:denominator-roots} in more
concrete terms.
Let $x_1,\dots,x_n$ be the elements of a
distinguished cluster~$\xx$ which correspond to the simple roots
$\alpha_1,\dots,\alpha_n$, respectively.

\begin{theorem}[\cite{ca2}]
\label{th:denominator-roots-restated}
For any root
$\alpha=c_1\alpha_1+\cdots+c_n\alpha_n \in\Phi_{\geq -1}\,$,
there is a unique cluster variable $x[\alpha]$ such that
\[
x[\alpha]=\frac{P_\alpha(x_1,\dots,x_n)}{x_1^{c_1}\cdots x_n^{c_n}}\,,
\]
where $P_\alpha$ is a polynomial in $x_1,\dots,x_n$ with nonzero
constant term; furthermore, any cluster variable is of this form.
\end{theorem}

In view of Theorems~\ref{th:denominator-roots} and
\ref{th:denominator-roots-restated}, the vertices
of the cluster complex $\Delta(\AA)$ (which we will also denote by
$\Delta(\Phi)$) can be identified with the almost positive roots in~$\Phi$.
Our goal is to describe the simplices of $\Delta(\Phi)$ explicitly
in root-theoretic terms.

\begin{theorem}[\cite{yga}]
\label{th:cluster-fan}
The~$n$ roots that label the
cluster variables in an arbitrary cluster
form a \hbox{$\ZZ$-basis} of the root lattice~$Q$.
The cones spanned by all such $n$-tuples of roots form a complete
simplicial fan 
in the ambient real vector space~$Q_\RR$.
\end{theorem}

Thus, the simplices of the cluster complex can be represented by
the cones of the simplicial fan in Theorem~\ref{th:cluster-fan}.
With some abuse of terminology, we denote this fan by the same
symbol $\Delta(\Phi)$.

\begin{figure}[ht]
\begin{center}
\setlength{\unitlength}{1.8pt}
\epsfbox{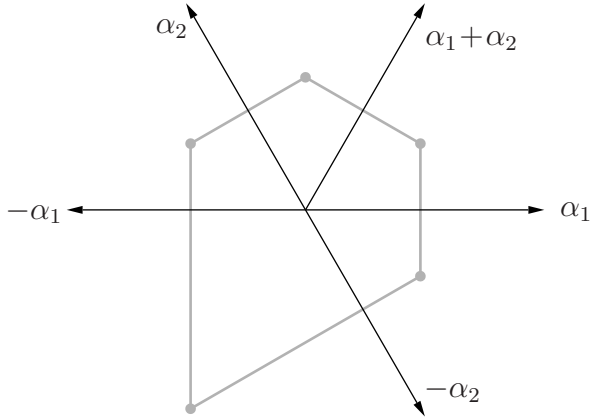}
\begin{picture}(0,0)(62.5,-47)
\thicklines
\put(57,0){\makebox(0,0){$\alpha_1$}}
\put(-57,0){\makebox(0,0){$-\alpha_1$}}
\put(35,36.5){\makebox(0,0){$\alpha_1\!+\!\alpha_2$}}
\put(-28,38.5){\makebox(0,0){$\alpha_2$}}
\put(31,-37.5){\makebox(0,0){$-\alpha_2$}}
\end{picture}
\end{center}
\caption{The fan $\Delta(\Phi)$ in type~$A_2$}
\label{fig:asso-A2}
\end{figure}

Figure~\ref{fig:asso-A2} shows the 
fan $\Delta(\Phi)$ for the special case of type~$A_2$.
This fan has $5$ maximal cones, spanned by the pairs
\begin{align*}
\{-\alpha_1,\alpha_2\},
\{\alpha_2,\alpha_1+\alpha_2\},
\{\alpha_1+\alpha_2,\alpha_1\},
\{\alpha_1,-\alpha_2\},
\{-\alpha_2,-\alpha_1\}
& \\
\quad \subset\Phi_{\geq -1}
=\{-\alpha_1,-\alpha_2,\alpha_1,\alpha_2,\alpha_1+\alpha_2\}. &
\end{align*}

The simplicial fan $\Delta(\Phi)$ in
Theorem~\ref{th:denominator-roots-restated}
is polytopal.

\begin{theorem}[\cite{cfz}]
\label{th:polytopality}
The fan $\Delta(\Phi)$ is the normal fan
of a simple \linebreak[3]
$n$-dimension\-al convex polytope, which is denoted
by~$P(\Phi)$ and called the generalized associahedron of type~$\Phi$.
\end{theorem}

Thus, the exchange graph of a cluster algebra of finite type
can be realized as the $1$-skeleton
of the generalized associahedron.


Turning to an explicit description of the cluster complex in
root-theoretic terms, we start with the following result from
\cite{yga}.

\begin{theorem}
\label{th:clique}
The cluster complex is a clique complex for its $1$-skeleton.
In other words, a subset $S\subset \Phi_{\geq-1}$ is a simplex in the
cluster complex if and only if every $2$-element subset of~$S$
is a \hbox{$1$-simplex}.
\end{theorem}

In order to describe the cluster complex, we therefore need only to
clarify which \emph{pairs} of almost positive roots label the edges of the cluster complex.
Without loss of generality, we assume that the root system $\Phi$
is irreducible; the general case can be obtained by taking direct
products.
Let $W$ be the Weyl group of~$\Phi$; it is a
Coxeter group generated by the simple reflections $s_1, \dots, s_n$
which correspond to the simple roots $\alpha_1,\dots,\alpha_n$.
Let $w_\circ$ be the element of maximal length in~$W$,
and $h$ the Coxeter number.
We also recall that $i \mapsto \varepsilon(i)$ is the sign function as in
Figure~\ref{fig:signs-dynkin}.




\begin{definition}\rm
\label{def:tau}
Define involutions $\tau_\pm:\Phi_{\geq
  -1}\to\Phi_{\geq -1}$ by
\end{definition}
\begin{equation*}
\label{eq:tau-pm-on-roots}
\tau_\varepsilon(\alpha) =
\begin{cases}
\displaystyle
\ \ \alpha & \text{if $\alpha = - \alpha_i\,$, with $\varepsilon(i)=- \varepsilon$;}
\\[.1in]
\Bigl(\displaystyle\prod_{\varepsilon(i) = \varepsilon} s_i\Bigr) \,(\alpha) & \text{otherwise.}
\end{cases}
\end{equation*}
For example, in type $A_2$
with $\varepsilon(1) = +1$ and $\varepsilon(2) = -1$, we get:
\begin{equation*}
\label{eq:A2-tau-tropical}
\hspace{-.1in}
\begin{array}{ccc}
-\alpha_1 & \stackrel{\textstyle\tau_+}{\longleftrightarrow}
~\alpha_1~ \stackrel{\textstyle\tau_-}{\longleftrightarrow}
~\alpha_1\,+\alpha_2~
\stackrel{\textstyle\tau_+}{\longleftrightarrow} ~\alpha_2~
\stackrel{\textstyle\tau_-}{\longleftrightarrow}
& -\alpha_2\,\, \\
\circlearrowright & & \circlearrowright \\ \tau_- & & \tau_+
\end{array}
\end{equation*}

\begin{theorem}[\cite{yga}]
\label{th:dihedral}
The order of $\tau_-  \tau_+$ is
$(h+2)/2$ if $w_\circ = -1$, and is $h+2$ otherwise.
Every $\langle \tau_-,\tau_+ \rangle$-orbit in  $\Phi_{\geq - 1}$ has a
nonempty intersection with $- \Pi$.
These intersections are precisely the $\langle - w_\circ \rangle$-orbits
in~$(- \Pi)$.
\end{theorem}

\pagebreak[3]

\begin{theorem}[\cite{yga,ca2}]
\label{th:compat}
There is a unique binary relation (called ``compatibility'') on
$\Phi_{\geq -1}$ that has the following two properties:
\begin{itemize}
\item
$\langle \tau_-,\tau_+ \rangle$-invariance:
if $\alpha$ and $\beta$ are compatible, then so are $\tau_\varepsilon \alpha$
and $\tau_\varepsilon\beta$ for $\varepsilon\in\{+,-\}$;
\item
a negative simple root $-\alpha_i$ is compatible with a root~$\beta$
if and only if the simple root expansion of~$\beta$ does not
involve~$\alpha_i$.
\end{itemize}
This compatibility relation is symmetric.
The clique complex for the compatibility relation is canonically
isomorphic to the cluster complex.
\end{theorem}

In other words (cf.\ Theorem~\ref{th:clique}), a subset of
$\Phi_{\geq -1}$ forms a simplex in the cluster complex
if and only if every pair of roots in this subset is compatible.

The above machinery allows us to explicitly describe a
generalized associahedron $P(\Phi)$ by a set of linear inequalities.

\begin{theorem}[\cite{cfz}]
\label{th:tau-inv-support-functions}
Suppose that a $(-w_\circ)$-invariant function \linebreak[3]
$F: -\Pi \to \RR$ satisfies the inequalities
\begin{equation*}
\label{eq:regular-dominant}
\sum_{i \in I} a_{ij} F(-\alpha_i) > 0 \quad \text{for all $j \in [1,n]$.}
\end{equation*}
Let us extend $F$ (uniquely) to a $\langle \tau_-,\tau_+
\rangle$-invariant function on~$\Phi_{\geq -1}\,$.
The generalized associahedron $P(\Phi)$ is then realized as the set
of points~$z$ in the dual space~$Q_\RR^*$ satisfying the linear inequalities
\begin{equation*}
\label{eq:inequalities-for-asso}
\langle z, \alpha \rangle\leq F(\alpha)\,,\ \text{for~all} \
\alpha\in\Phi_{\geq -1}\, .
\end{equation*}
\end{theorem}

As an example of a function~$F$ satisfying the conditions in
Theorem~\ref{th:tau-inv-support-functions},
one can take the function whose value
$F(-\alpha_i)$ is equal to the coefficient of a simple coroot~$\alpha_i^\vee$
  in the half-sum 
of all positive coroots.

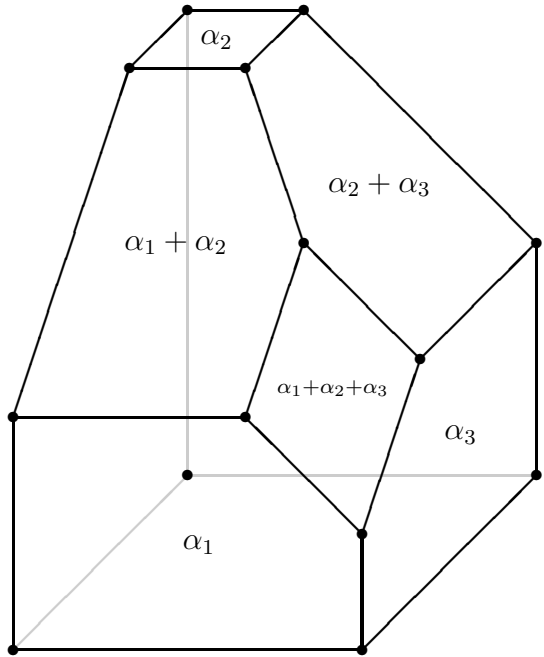
\begin{figure}[ht!]
\begin{center}
\setlength{\unitlength}{2.2pt}
\begin{picture}(90,110)(0,0)
\thicklines

\verylight{
\put(30,30){\line(0,1){80}}
\put(30,30){\line(1,0){60}}
\put(0,0){\line(1,1){30}}
 }
\put(0,0){\line(1,0){60}}
\put(0,0){\line(0,1){40}}
\put(60,0){\line(0,1){20}}
\put(60,0){\line(1,1){30}}
\put(0,40){\line(1,0){40}}
\put(0,40){\line(1,3){20}}
\put(60,20){\line(-1,1){20}}
\put(60,20){\line(1,3){10}}
\put(90,30){\line(0,1){40}}
\put(40,40){\line(1,3){10}}
\put(70,50){\line(1,1){20}}
\put(70,50){\line(-1,1){20}}
\put(50,70){\line(-1,3){10}}
\put(90,70){\line(-1,1){40}}
\put(20,100){\line(1,0){20}}
\put(20,100){\line(1,1){10}}
\put(30,110){\line(1,0){20}}
\put(40,100){\line(1,1){10}}

\put(35,105){\makebox(0,0){$\alpha_2$}}
\put(28,70){\makebox(0,0){$\alpha_1+\alpha_2$}}
\put(63,80){\makebox(0,0){$\alpha_2+\alpha_3$}}
\put(55,45){\makebox(0,0){$\scriptstyle \alpha_1+\alpha_2+\alpha_3$}}
\put(32,18){\makebox(0,0){$\alpha_1$}}
\put(77,37){\makebox(0,0){$\alpha_3$}}

\put(0,0){\circle*{2}}
\put(60,0){\circle*{2}}
\put(60,20){\circle*{2}}
\put(30,30){\circle*{2}}
\put(90,30){\circle*{2}}
\put(0,40){\circle*{2}}
\put(40,40){\circle*{2}}
\put(70,50){\circle*{2}}
\put(50,70){\circle*{2}}
\put(90,70){\circle*{2}}
\put(20,100){\circle*{2}}
\put(40,100){\circle*{2}}
\put(30,110){\circle*{2}}
\put(50,110){\circle*{2}}

\end{picture}
\end{center}
\caption{
The type $A_3$ associahedron}
\label{fig:A3asso-poly}
\end{figure}

\begin{example}{\rm
In type~$A_3$, Theorem~\ref{th:tau-inv-support-functions} is illustrated in
Figure~\ref{fig:A3asso-poly}, which shows a $3$-dimensional simple
polytope
given by the inequalities
\[
\begin{array}{rcl}
\max(-z_1\,,\,
-z_3\,,\,
z_1\,,\,
z_3\,,\,
z_1+z_2\,,\,
z_2+z_3
)&\!\!\leq\!\!& 3/2\,,\\[.1in]
\max(-z_2\,,\,
z_2\,,\,
z_1+z_2+z_3
)&\!\!\leq\!\!& 2\,.
\end{array}
\]
}
\end{example}

For the classical types, the clusters and generalized assohiahedra
have concrete combinatorial realizations in terms of triangulations of
regular polygons.
This realization is discussed
in detail in \cite{yga, cfz} and further explored in~\cite{ca4}.
Here we restrict ourselves to an outline of the type~$A_n$
construction,
where we recover the classical associahedron, originally introduced by
Stasheff~\cite{stasheff}.
More precisely, the cluster complex of type $A_n$ is isomorphic
to the dual complex of an $n$-dimensional associahedron, as explained
in Example~\ref{example:compat-An} below.

\begin{example}
\label{example:compat-An}
The cluster complex of type~$A_n$
can be naturally identified with the simplicial complex
whose simplices are all collections of non-crossing diagonals of a
regular $(n+3)$-gon~$\mathbf{P}_{n+3}\,$ (cf.\ Section~\ref{sec:triangulations}).
Thus, we need to identify
the roots in $\Phi_{\geq -1}$ with the diagonals of $\mathbf{P}_{n+3}\,$.
This is done as follows.
The roots in $- \Pi$ correspond to the diagonals on the ``snake''
shown in Figure~\ref{fig:octagon-snake}.

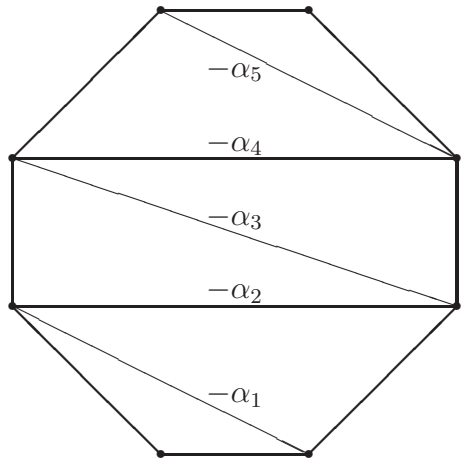
\begin{figure}[ht]
\begin{center}
\setlength{\unitlength}{2.8pt}
\begin{picture}(60,60)(0,0)
\thicklines
  \multiput(0,20)(60,0){2}{\line(0,1){20}}
  \multiput(20,0)(0,60){2}{\line(1,0){20}}
  \multiput(0,40)(40,-40){2}{\line(1,1){20}}
  \multiput(20,0)(40,40){2}{\line(-1,1){20}}

  \multiput(20,0)(20,0){2}{\circle*{1}}
  \multiput(20,60)(20,0){2}{\circle*{1}}
  \multiput(0,20)(0,20){2}{\circle*{1}}
  \multiput(60,20)(0,20){2}{\circle*{1}}

\thinlines \put(0,20){\line(1,0){60}} \put(0,40){\line(1,0){60}}
\put(0,20){\line(2,-1){40}} \put(0,40){\line(3,-1){60}}
\put(20,60){\line(2,-1){40}}

\put(30,8){\makebox(0,0){$-\alpha_1$}}
\put(30,22){\makebox(0,0){$-\alpha_2$}}
\put(30,32){\makebox(0,0){$-\alpha_3$}}
\put(30,42){\makebox(0,0){$-\alpha_4$}}
\put(30,52){\makebox(0,0){$-\alpha_5$}}


\end{picture}
\end{center}
\caption{The ``snake'' in type $A_5$}
\label{fig:octagon-snake}
\end{figure}

Each positive root $\alpha_i + \alpha_{i+1} + \cdots + \alpha_j$
corresponds to the
unique diagonal that crosses precisely the diagonals $- \alpha_i,
- \alpha_{i+1}, \ldots,- \alpha_j$ from the snake.
As shown in \cite{yga}, two roots are compatible if and only if the
corresponding diagonals do not cross each other at an interior point.
Consequently, an $n$-tuple of roots
$\beta_1,\dots,\beta_n \in \Phi_{\geq -1}$
forms a maximal simplex in the cluster complex
if and only if the corresponding diagonals do not cross each
other---i.e., they form a triangulation of~$\mathbf{P}_{n+3}\,$.
(In particular, in agreement with Example~\ref{ex:type-An},
the ``snake" triangulation is associated 
with a distinguished cluster.)
As a result, the exchange graph for a cluster algebra of
type~$A_n$ becomes identified with the graph that already appeared
in Section~\ref{sec:triangulations} (with triangulations
of~$\mathbf{P}_{n+3}$ as vertices, and diagonal flips as edges).
The type~$A_3$ case is illustrated by Figures~\ref{fig:A3assoc_poly}
and~\ref{fig:A3asso-poly}.
\end{example}

\begin{figure}[ht]
\centerline{\scalebox{0.9}{\epsfbox{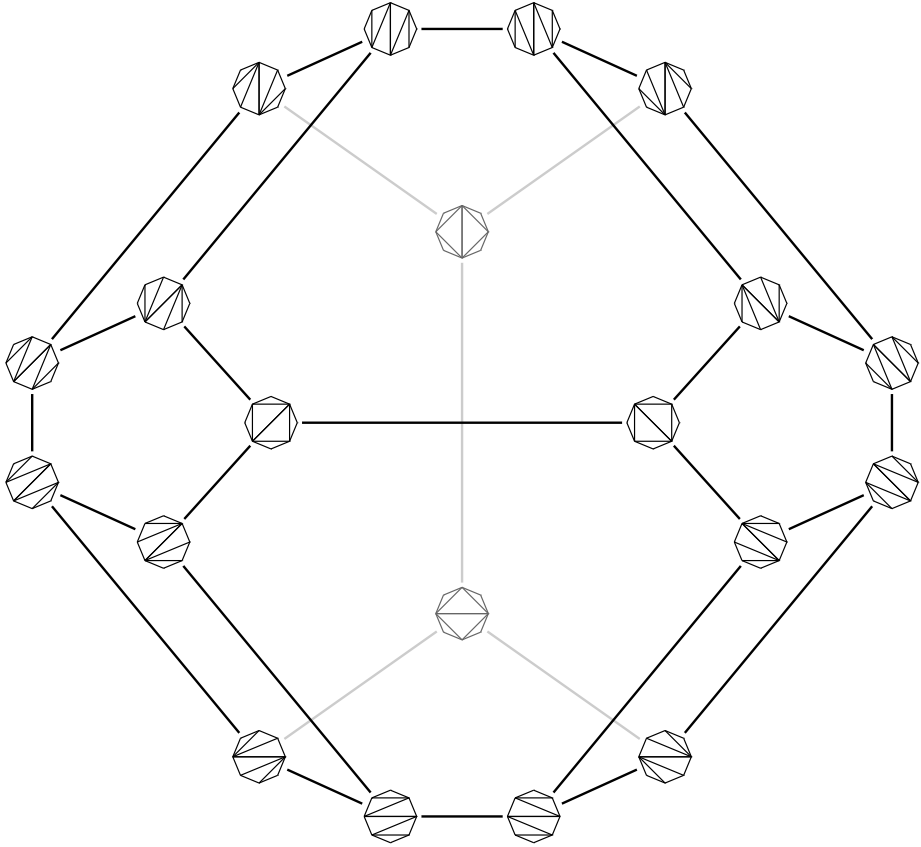}}}
\caption{The $3$-dimensional cyclohedron}
\label{B3assoc}
\end{figure}

In types~$B_n$ and~$C_n$, the cluster complex is isomorphic to the
dual complex of the $n$-dimensional \emph{cyclohedron}, or
Bott-Taubes polytope~\cite{bott-taubes}.
The vertices of this polytope are labeled
by the centrally symmetric triangulations of a regular
$(2n+2)$-gon~$\mathbf{P}_{2n+2}\,$,
while each edge corresponds to a pair of centrally symmetric flips,
or to a single flip inside a rectangle.
(We omit the description of the dictionary between the almost positive
roots in a root system of type $B_n$ or $C_n$ and the (pairs of) diagonals
in~$\mathbf{P}_{2n+2}\,$, referring the reader to~\cite{ca2} or~\cite{yga}.)
See Figure~\ref{B3assoc}.

Similar description can be provided for the generalized associahedron
of type~$D_n$; see~\cite{ca2,yga}.

These combinatorial models for the cluster complexes of classical
types $ABCD$ make the study of the corresponding cluster algebras much more
tractable.
In particular, these models are utilized in the proof of
Theorem~\ref{th:conjectures-classical}.
It is therefore very desirable to find a solution of the following
tantalizing problem.

\begin{problem}
\label{prob:concrete-exceptional}
Construct combinatorial models for the generalized
associahedra of exceptional types.
\end{problem}

\subsection{Numerology}

As discussed in Example~\ref{example:compat-An},
the clusters in a cluster algebra of
type~$A_n$ are in bijection with the triangulations of a convex
$(n+3)$-gon by $n$ non-crossing diagonals.
The number of such triangulations is well known to be the Catalan
number $\frac{1}{n+2}\binom{2n+2}{n+1}$.
The following theorem gives a general formula for the number of
clusters in a cluster algebra of arbitrary finite type.

\begin{theorem}[\cite{yga}]
\label{th:N-thru-exponents}
The number of clusters in a cluster algebra of finite type
associated with a root system~$\Phi$
(equivalently, the number of vertices of the corresponding
generalized associahedron) is equal to
\begin{equation}
\label{eq:N-thru-exponents}
N(\Phi) 
{=} \prod_{i=1}^n \frac{e_i + h + 1}{e_i +
  1} \,,
\end{equation}
where $e_1, \dots, e_n$ are the exponents of~$\Phi$,
and $h$ is the Coxeter number.
\end{theorem}

Figure~\ref{W-Cat} shows the values of $N(\Phi)$ for all
irreducible crystallographic root systems~$\Phi$.

\begin{figure}[ht]
\begin{center}
\begin{tabular}{|c|c|c|c|c|c|c|c|c|}
\hline
$A_n$ &
$B_n,C_n
$ & $D_n$ & $E_6$ & $E_7$ & $E_8$ & $F_4$  & $G_2$ \\
\hline
&&&&&&&\\[-.1in]
$\textstyle\frac{1}{n+2} \binom{2n+2}{n+1}$&
$\textstyle\binom{2n}{n}$ 
&   $\textstyle\frac{3n-2}{n}\binom{2n-2}{n-1}$ 
&
$\textstyle 833$
&
$\textstyle 4160$
&
$\textstyle 25080$
&
$\textstyle 105$
& $\textstyle 8$ \\[.05in]
\hline
\end{tabular}
\end{center}
\caption{The numbers $N(\Phi)$ 
}
\label{W-Cat}
\end{figure}

The numbers $N(\Phi)$ given by \eqref{eq:N-thru-exponents}
can be thought of as generalizations of the
Catalan numbers to an arbitrary Cartan-Killing type.
These numbers are known to count
\cite{bessis, cellini-papi, djokovic, haiman, picantin,
  Reiner-noncrossing, shi} a variety of
combinatorial objects related to the root system~$\Phi$.
In particular, $N(\Phi)$ appears as the number of:
\begin{itemize}
\item
\emph{$ad$-nilpotent ideals} in a semisimple Lie algebra;
\item
antichains in the \emph{root poset};
\item
positive regions of the \emph{Shi arrangement};
\item
$W$-orbits in the quotient $Q^\vee/(h+1)Q^\vee$ of the dual root
lattice;
\item
conjugacy classes of elements of a semisimple Lie
group which have a finite order that divides~$h+1$;
\item
\emph{non-crossing partitions} of appropriate type.
\end{itemize}
We refer the reader to \cite[Lecture~5]{pcmi} for a
discussion of these enumerative formulas and their $q$-analogues,
as well as for bibliographic directions.

As noted in \cite[Section~5.3]{pcmi},
substantial part of the combinatorial theory of
generalized associahedra survives in the case of non-cryctallographic
finite root systems $I_2(n)$, $H_3$, and~$H_4$,
including the construction of the clique complex in
Theorem~\ref{th:compat} and formula~\eqref{eq:N-thru-exponents}
for the number of root clusters---but excluding cluster algebras as
such.

\subsection{Concluding example}

We close this survey with an example of a cluster
algebra of finite type that arises from the geometric construction
of Section~\ref{sec:tp-dbc}.

\begin{example}
\label{ex:sl3-D4}
[Coordinate ring of a double Bruhat cell
$G^{w_\circ,w_\circ}\subset SL_3$]
Recall that the open double Bruhat cell $G^{w_\circ,w_\circ}\subset SL_3(\CC)$
consists of all complex $3\times 3$ matrices $x=(x_{ij})$
of determinant~$1$ whose minors
\begin{equation}
\label{eq:sl3-w0w0-nonvanishing}
x_{13},
{\ \ }
\Bigl|\!\Bigl|\begin{array}{cc}
x_{12} & x_{13} \\
x_{22} & x_{23}
\end{array}\Bigr|\!\Bigr|,
{\ \ }
x_{31},
{\ \ }
\Bigl|\!\Bigl|\begin{array}{cc}
x_{21} & x_{22} \\
x_{31} & x_{32}
\end{array}\Bigr|\!\Bigr|
\end{equation}
are nonzero.
The coordinate ring $\QQ[G^{w_\circ,w_\circ}]$ is a cluster algebra of
type~$D_4$
over the ground ring generated by the minors in~\eqref{eq:sl3-w0w0-nonvanishing} and
their inverses.
It has $16$~cluster 
variables, corresponding to the $16$ roots in~$\Phi_{\geq -1}$:

\begin{itemize}
\item[(i)]
$14$ (among the $19$ total) minors of~$x$, namely,
all excluding $\det(x)$ and the $4$~minors
in~\eqref{eq:sl3-w0w0-nonvanishing};

\item[(ii)]
$x_{12}x_{21}x_{33}-x_{12}x_{23}x_{31}-x_{13}x_{21}x_{32}+x_{13}x_{22}x_{31}$;

\item[(iii)]
$x_{11}x_{23}x_{32}-x_{12}x_{23}x_{31}-x_{13}x_{21}x_{32}+x_{13}x_{22}x_{31}$.
\end{itemize}

They form $50$~clusters of size~$4$, one for each of
the $50$ vertices of the type~$D_4$ associahedron.

These clusters give rise to 50 different \emph{total positivity criteria}:
a matrix
$x\in SL_3$ is totally positive if and~only if the following
$8$~functions are positive at~$x$:
the $4$~elements of any given cluster and the $4$~minors in~\eqref{eq:sl3-w0w0-nonvanishing}.
\end{example}

\end{document}